\documentclass[11pt]{amsart}
\usepackage{palatino}
\usepackage[all]{xy,xypic}
\usepackage{amsfonts,amsmath,amssymb,amscd,enumerate}
\usepackage{tikz-cd}
\usepackage{oldgerm}
\usepackage[a4paper,left=12em,right=12em,top=7em,bottom=7em]{geometry}
\usepackage{hyperref}

\makeatletter
\@namedef{subjclassname@2020}{%
	\textup{2020} Mathematics Subject Classification}
\makeatother

\newcommand{\ra}{\rightarrow}		
\newcommand{\lra}{\longrightarrow}
\newcommand{\by}[1]{\stackrel{#1}{\ra}}

\newcommand{\remove}[1]{}

\newcommand{\ol}{\overline}		

\newcommand{\iso}{\by \sim}

\newtheorem{theorem}{Theorem}[section]
\newtheorem{proposition}[theorem]{Proposition}
\newtheorem{lemma}[theorem]{Lemma}
\newtheorem{definition}[theorem]{Definition}
\newtheorem{corollary}[theorem]{Corollary}
\newtheorem{conjecture}[theorem]{Conjecture}
\newtheorem{example}[theorem]{Example}

	\newcommand{\bz}{\mbox{$\mathbb Z$}}

\newcommand{\MI}{\mbox{$\mathfrak I$}}

\newcommand{\ms}{\mbox{$\mathfrak s$}}

\newcommand{\ot}{\mbox{\,$\otimes$\,}}	
\newcommand{\op}{\mbox{$\oplus$}}
\newcommand{\Spec}{\text{Spec}}	
\newcommand{\hh}{\text{ht}}
\newcommand{\rank}{\text{rank}}

\newcommand{\Aut}{\mbox{\rm Aut\,}}

\newcommand{\Um}{\text{Um}}		
\newcommand{\SL}{\text{SL}}
\newcommand{\GL}{\text{GL}}		
	
\newcommand{\sur}{\twoheadrightarrow}
\newcommand{\bp}{\begin{proposition}}
	\newcommand{\ep}{\end{proposition}}
\newcommand{\bl}{\begin{lemma}}
	\newcommand{\el}{\end{lemma}}
\newcommand{\bt}{\begin{theorem}}
	\newcommand{\et}{\end{theorem}}
\newcommand{\bc}{\begin{corollary}}
	\newcommand{\ec}{\end{corollary}}
\newcommand{\bd}{\begin{definition}}
	\newcommand{\ed}{\end{definition}}
\newcommand{\bco}{\begin{conjecture}}
	\newcommand{\eco}{\end{conjecture}}
\newcommand{\bma}{\begin{bmatrix}}
	\newcommand{\ema}{\end{bmatrix}}

\renewcommand{\k}{\overline{\mathbb F}_p}

\def\rmk{\refstepcounter{theorem}\paragraph{{\bf Remark} \thetheorem}}
\def\proof{\paragraph{Proof}}

\def\notation{\paragraph{\bf Notation}}
\def\quest{\refstepcounter{theorem}\paragraph{{\bf Question} \thetheorem}}

\title[Splitting criteria for projective modules over polynomial algebras]{Splitting criteria for projective modules\\
	over polynomial algebras}
\date{\today}
\author{Sourjya Banerjee and Mrinal Kanti Das}
\address{(Sourjya Banerjee) Stat-Math Unit, Indian Statistical Institute, 203, B.T. Road, Kolkata-700108, West Bengal, INDIA}
\curraddr{(Sourjya Banerjee) The Institute of Mathematical Sciences, HBNI, C.I.T. Campus, Tharamani, Chennai  600113, India}
	\email{sourjyab@imsc.res.in, sourjya91@gmail.com}

\address{(Mrinal Kanti Das) Stat-Math Unit, Indian Statistical Institute, 203, B.T. Road, Kolkata-700108, West Bengal, INDIA} 
\email{mrinal@isical.ac.in, das.doublelife@gmail.com}

\begin{document}
	\maketitle

	\begin{abstract}	
	This article investigates the splitting problem for finitely generated projective modules $P$ over affine algebras over algebraically closed fields and their polynomial extensions. We then address an open question due to M. Roitman on monic inversion principle for projective modules and prove it in the affirmative  for finitely generated rings. For affine algebras over $\ol{\mathbb{F}}_p$, we prove a monic inversion principle for ideals. We also exhibit some applications.
	\end{abstract}

	\section{Introduction}
The purpose of this article is to understand the splitting behavior of a finitely generated projective module (of a constant rank) over a polynomial algebra in terms of its generic sections.  Our approach varies in two distinct yet related themes, which we present one by one.

{\bf All the rings below are assumed to be commutative and Noetherian. The modules are assumed to be finitely generated}.

\subsection{Splitting of projective modules}
We first recall a remarkable result due to N. Mohan Kumar. Although the result is initially stated  differently in \cite[Theorem 1]{NMK84}, a closer look at his proof reveals that he essentially establishes the following theorem. The notation $\mu(-)$ stands for the minimal number of generators.

\bt\cite{NMK84}\label{NMKSC}
Let $A$ be a ring of dimension $d \ge 2$ and $P$ be a projective $A$-module of rank $d$. Let $I$ be an ideal of $A$ and   $\alpha:P\sur I$ be a surjection. 
\begin{enumerate}
	\item
	If  $A$ is a reduced affine algebra  over an algebraically closed field $k$, then
	$\mu(I)=d$ implies that $P\simeq Q\op A$ (for some $A$-module $Q$).
	\item
	$P\simeq Q\op A$ and $\det(P)$ trivial imply that $\mu(I)=d$ (for any ring $A$).
\end{enumerate}
\et 

The first assertion played a key role in M. P. Murthy's seminal work \cite{Mu}, where he proved that the vanishing of the top Chern class $c_d(P)$ in $\text{CH}^d(\Spec(A))$ dictates the splitting behavior for $P$, when $A$ is a smooth affine algebra over an algebraically closed field $k$. The smoothness assumption has subsequently been removed by A. Krishna in \cite{Kr}. 

%The second assertion remains valid for $P$ with a nontrivial determinant, thanks to the following: \cite[Theorem 3.7]{Mu} and \cite[Corollary 1.5]{Kr}.

Therefore, it is natural  to inquire whether one can establish a criterion (similar to (1) above) for the splitting of $P$ when $\rank(P) < \dim(A)$. Unfortunately, progress in addressing even the scenario where $\rank(P) = \dim(A) - 1$, has been notably slow. Recently, A. Asok and J. Fasel \cite{AF1, AF2} tackled the case $\rank(P) = \dim(A) - 1$, where $A$ is a smooth affine algebra of dimensions $3$ and $4$ over an algebraically closed field of characteristic unequal to $2$. However, their assertion involves  the vanishing of  $c_{d-1}(P)$ in $\text{CH}^{d-1}(\Spec(A))$. However, the analogue of (1) is not apparent from their results and we believe that it would be daunting task to relate the two. Even then, their methods are	quite involved and they have to employ some sophisticated machinery to prove these results, which gives an idea about the depth of the problem. 

%In this context, we investigate the scenario $\rank(P) = \dim(A) - 1$, considering $A$ as a polynomial algebra over an algebraically closed field of characteristic unequal to $2$.  We prove the following theorem [Theorem \ref{BRQRLT} and Corollary \ref{EUESFM}].

%	 and if $P$ is a projective $A$-module of rank $d-1$, then $c_{d−1}(P)$ vanishes in $\text{CH}_{d−1}(X)$ if and only if . Their methods are
%	quite involved and they have to employ some sophisticated machinery to prove these
%	results, which gives an idea about the depth of the problem.	\cite{AF1} \cite{AF2} In this regard, the first and obvious set up to explore is the polynomial algebra. We prove the following theorem. 

In a spirit similar to Mohan Kumar's result, we prove the following theorem [Theorem \ref{BRQRLT} and Corollary \ref{EUESFM}].

\bt\label{MT1I}
Let $A$ be an affine algebra of dimension $d \ge 2$ over an algebraically closed field $k$, where $\text{char}(k)\not=2$.
Let $P$ be a projective $A[T]$-module of rank $d$. Let $I$ be an ideal of $A[T]$ and  $\alpha:P\sur I$ be a surjction. If 
$\mu(I)=d$, then $P\simeq Q\op A[T]$ (for some $A[T]$-module $Q$). Consequently, if $P$ is stably free then $P$ always splits off a free summand of rank one. 	
\et

It is worth pointing out that it is in no way a simple translation of Mohan Kumar's methods to the above set up. Prior to us,
 S. M. Bhatwadekar and R. Sridharan proved this result \cite[Theorem 4.5]{BRS01} when $\text{char}(k)=0$ and the determinant of $P$ is trivial and the proof is quite involved. One may observe that the same proof can be extended in the case when $\text{gcd}(\text{char}(k), d!)=1$. For us, to extend the hypotheses to the case $\text{char}(k) \neq 2$ and an arbitrary determinant for $P$, a completely different line of reasoning was required [cf. Section \ref{5}]. On the other hand, compared to \cite{AF1, AF2}, our method employs simple algebraic tools and delves into the polynomial structure in depth. When $k=\k$,  we show in Theorem \ref{MKLIPR} that the implication in the other direction also holds. We prove:

\bt\label{fp}
Let $A$ be an affine algebra of dimension $d \ge 2$ over $\k$ (no restriction on $p$).
Let $P$ be a projective $A[T]$-module of rank $d$ with trivial determinant. Let 
$I$ be an ideal in $A[T]$ and  $\alpha:P\sur I$ be a surjection. 
If $P\simeq Q\op A[T]$ (for some $A[T]$-module $Q$),
then $\mu(I)=d$. 
\et

%Theorem \ref{fp} is also true if we replace $\k$ by  any infinite field and assume that $(d-1)!$ is invertible in $A$.
%We indicate the method in Remark \ref{infinitefield}.

\smallskip

In \cite{MKDIMRN}, M. K. Das improved N. Mohan Kumar's result for the case when $k=\k$, under the much stronger assumptions that either $\frac{1}{(d-1)!}\in A$ or $A$ is smooth. We further refine the improvement presented in \cite{MKDIMRN} in the following manner (for details we refer to \ref{aafpbar}).

\bt
Let $p\not= 2,3$ and $A$ be an affine algebra of dimension $d\ge 5$ over $\k$. Let $P$
be a projective $A$-module of rank $d-1$ with trivial determinant. Let $\alpha: P\sur I$ be a surjection. Then $P$ splits off a free summand of rank one if and only if $\mu(I)=d-1$.  
\et

\subsection{Monic inversion principle}
 M. Roitman \cite{QOMR} posed the following splitting problem for projective modules over polynomial algebras.\\

\quest \label{question3}  Let $A$ be a commutative Noetherian ring of dimension $d\ge 2$ and $P$ be a projective $A[T]$-module of rank $d$.  Suppose that there exists a surjection $\alpha:P\sur I$, where $I\subset A[T]$ is an ideal. Moreover assume that $I$ contains a monic polynomial in $T$. Then does $P$ split a free summand of rank one?\\

The above question remains open in general. In the same paper \cite{QOMR}  M. Roitman answered the question affirmatively when $A$ is local. Later S. M. Bhatwadekar and R. Sridharan \cite[Theorem 3.4]{BRS01} gave an affirmative answer when the ring contains an infinite field. In exploring the remaining cases, we provide an affirmative answer to the above question over finitely generated $\bz$-algebras in the following manner.

\bt(Theorem \ref{RQOZ})
Let $A$ be a finitely generated $\mathbb{Z}$-algebra of dimension $d\ge 2$ such that there exists an integer $n\ge 2$  that is invertible in $A$.  Let $P$ be a projective $A[T]$-module of rank $d$ with trivial determinant,  and $I\subset A[T] $ be an ideal of height $d$ containing a monic polynomial in $T$. Suppose that there exists a surjection $\alpha:P\sur I$. Then $P$ splits off a free summand of rank one.
\et

On a related note, let us recall  the following question on \emph{``monic inversion principle"} for ideals, posed by S. M. Bhatwadekar. Here $A(T)$ is the ring obtained from $A[T]$ by inverting all the monic ploynomials.

\smallskip

\quest \label{question1} Let $A$ be a commutative Noetherian ring of dimension $d\ge 2$ and $I\subset A[T]$ be an ideal such that 
$\hh(I)=\mu(I/I^2)=d$. Moreover assume that $I=\langle f_1,...,f_d\rangle  +I^2$. Suppose that there exists $F_i\in IA(T)$ such that $IA(T)=\langle F_1,...,F_d\rangle  $, with $F_i-f_i\in I^2A(T)$. Then does there exist $g_i\in I$, such that $I=\langle g_1,...,g_d\rangle  $, with $g_i-f_i\in I^2$?

\medskip

Question \ref{question1} has a negative answer for $d=2$ ( see \cite[Example 3.15]{BR} ). Whenever $d\ge 3$ this question has an affirmative answer in the following cases :
\begin{itemize}
	\item $A$ is a local ring \cite[ Proposition 5.8(1)]{MKD1}.
	\item $A$ is an affine domain over an algebraically closed field of characteristic zero \cite[Proposition 5.8(2)]{MKD1}.
	\item $A$ is a regular domain which is essentially of finite type over an infinite perfect field $k$ of characteristic unequal to $2$ \cite[Theorem 5.11]{DTZ2}.
	\item $A$ is a real affine algebra such that: either there are no real maximal ideals, or the intersection of all real maximal ideals has height at least one \cite{BSOU}.
	
\end{itemize}

We prove the following result, using methods  entirely different from the works mentioned above.

\bt\label{II}(Theorem \ref{INJECG})
Let $R$ be an affine algebra over $\ol{\mathbb{F}}_p$ of dimension $d\ge 2$ and $I\subset R[T]$ be an ideal such that $\hh(I)=\mu(I/I^2)=d$. Moreover assume that $I=\langle f_1,...,f_d\rangle  +I^2$. Suppose that there exists $F_i\in IR(T)$ such that $IR(T)=\langle F_1,...,F_d\rangle  $, with $F_i-f_i\in I^2R(T)$. Then there exists $g_i\in I$, such that $I=\langle g_1,...,g_d\rangle  $, where $g_i-f_i\in I^2$.
\et

%	As an application of the above result we prove some addition and subtraction principles on $R[T]$, where $R$ is an affine algebra over $\k$ of dimension $d\ge 2$. Let $P$ be a projective $R[T]$-module of rank $d$ with a trivial determinant and $\chi:R[T]\cong \wedge^dP$ be an isomorphism. We then assign a \textit{``local orientation"} $(I,\omega_I)\in E^d(R[T])$ to the pair $(P,\chi)$ and show that $P$ splits a free summand of rank one if and only if the image of $(I,\omega_I)$ in the group $E^d(R[T])$ vanishes. Here $E^d(R[T])$ is the $d$-th Euler class group of $R[T]$. Moreover assuming $(d-1)!$ is invertible, we show that $E^d(R[T])$ is the obstruction group to dictate the splitting problem of $P$. In \cite{MKDIMRN}, M. K. Das showed that $E^{d-1}(A)$ is the precise obstruction group to dictate the splitting problem, where $A$ is a $d$-dimensional \emph{\bf smooth} $\k$-algebra. Here we remark that, over smooth polynomial extensions since projective modules are extended, removing the ``smoothness" hypothesis is crucial.   

\section*{Acknowledgment}
This article was a part of the first named author's Ph.D. thesis. He would like to express his gratitude to the Stat-Math Unit of the Indian Statistical Institute, Kolkata, for their invaluable support.

	\section{Preliminaries}
	In this section, we recollect various results from the literature.  We start with a lemma due to N. Mohan Kumar, recast slightly to suit our needs.
	
	\bl\label{MKL}\cite{NMKL}
	Let $A$ be a ring and $I$ be an ideal of $A$. Let $J, K$ be
	ideals of $A$ contained in $I$ such that $K\subset I^2$ and $I=J+K$. Then there exists $e\in K$ such that $e(1-e)\in J$ and $I=\langle J,e\rangle  .$
	\el

	The following lemma is known as ``Moving Lemma". We restate it to suit our requirements. The proof is given in \cite[Corollary 2.14]{MPHIL}, (one just needs to use \cite[Theorem 2.4]{MKD} in the appropriate places to establish the second assertion). Hence we skip the proof to avoid repeating the same arguments.

	\bl(Moving Lemma)\label{ML}
	Let $A$ be a  ring of dimension $n\geq 2$ and $I\subset A$ be an ideal such that $\mu(I/I^2)=n=\hh(I)$. Let $I=\langle a_1,...,a_n\rangle  +I^2$. Then there exists an ideal $J\subset A$, either of height $n$ or $J=A$, with the property that $I\cap J=\langle b_1,...,b_n\rangle   $, with $a_i-b_i\in I^2$ and $I+J=A$. Moreover, if $A$ is an affine algebra over $\k$ of dimension $n+1$ and $\hh(I)=n$, then $J$ can be chosen to be co-maximal with any ideal of height $\ge 1$.
	\el
	%\proof Using Lemma \ref{MKL} get $e\in I^2$ be such that $I=\langle a_1,...,a_d,e\rangle  $ with $e(1-e)\in \langle a_1,...,a_n\rangle  $. By (\cite{SMBB3}, Corollary 2.13) there exists $\lambda_i\in A$ such that $\hh(\langle b_1,...,b_n\rangle  )_e\ge n$, where $b_i=a_i+e\lambda_i$, for $i=1,...,n$. Let $J=\langle b_1,...,b_n,1-e\rangle  $
	The next theorem is an accumulation of results due to different authors. For a proof one can see \cite[Corollary 2.4]{MKDIMRN}.
	\bt\label{MKDLE}
	Let $R$ be an affine algebra over $\k$, and let $I \subset R$ be an ideal such that
	$\dim(R/I) \le 1$. Then, we have the following assertions:\\
	$(1)$ The canonical map $\SL_n(R)\sur SL_n(R/I)$ is surjective for $n \ge 3$.\\
	$(2)$ If $\dim(R) = 3$, then the canonical map $\SL_2(R) \sur \SL_2(R/I)$ is surjective.
	\et
	
	\smallskip
	
	\bd
	Let $R$ be a ring. A vector $(a_1,...,a_n)\in R^n$ is called a unimodular row of length $n$ if there exists $(b_1,...,b_n)\in R^n$ such that $a_1b_1+...+a_nb_n=1$. We will denote $\Um_n(R)$ by the set of all unimodular rows of length $n$.
	\ed

	Then next lemma is frequently used in this article. It can be deduced following  the proof of \cite[Lemma 5.3]{SMBB3}. For the sake of completeness we give a proof. Before that we recall a notation which will be used frequently.
	
	\notation Let $R$ be a ring. For any two square matrices $M\in M_{m}(R)$ and $N\in M_{n}(R)$, by $M\perp N$ we denote the following matrix. $$\begin{pmatrix}
		M & {0} \\
		0 & N
	\end{pmatrix}\in M_{{m+n}}(R)$$

	\bl\label{NLTP}
	Let $R$ be an affine algebra over $\ol{\mathbb{F}}_p$ of dimension $d+1\ge 3$ and $I\subset R$ be an ideal such that $\hh(I)=\mu(I/I^2)=d$. Let $f\in R$ be a unit modulo $I$. Assume that, there is a surjection $\omega_I:(R/I)^d\sur I/I^2$ which induces $I=\langle f_1,...,f_d\rangle  +I^2$. Let $\alpha\in \GL_d(R/I)$, such that $\det(\alpha)=\widetilde{f^2}$, where `tilde' denotes going modulo $I$. Suppose that, the map $\omega_I\alpha:(R/I)^d\sur I/I^2$ induces $I=\langle g_1,...,g_d\rangle  +I^2$. If $ f_1,...,f_d$ can be  lifted to a set of $d$-generators of $I$, then $ g_1,...,g_d$ can also be lifted to a set of $d$-generators of $I$.
	\el
	\proof  Let $I=\langle h_1,...,h_d\rangle  $ where $h_i-f_i\in I^2$. It follows from prime avoidance lemma (cf. \cite[Proposition 7.1.2]{IR}) that, after an elementary transformation, we may assume $\hh(\langle h_1, \dots, h_i \rangle) = i$ for each $i = 1, \dots, d$. Let $B=R/\langle h_3,...,h_d\rangle  $ and `bar' denote going modulo $\langle h_3,...,h_d\rangle  $. Then $\dim(B)\le 3$.
	\par Since $f$ is a unit modulo $I$, we have $g\in R$ such that $fg-1\in I$. Note that $(\ol g^2,\ol h_2,-\ol h_1)\in \Um_3(B)$. By a result Swan-Towber \cite{SWANT} the unimodular row $(\ol g^2,\ol h_2,-\ol h_1)$ is completable to an invertible matrix in $\SL_3(B)$. Using  \cite[Lemma 5.2]{SMBB3} we get $\tau\in M_2(B)$ such that $(\ol h_1,\ol h_2)\tau =(\ol h_1',\ol h_2')$, where $\ol I= \langle \ol h_1',\ol h_2'\rangle  $ and $\det(\tau')-f^2\in I$, for some $\tau'\in M_2(R)$ which is a lift of $\tau$. With an abuse of notation from now onward throughout the article, we shall not distinguish between $\tau$ and $\tau'$, and will only denote $\det(\tau)-f^2\in I$ to convey the same.

	\par Thus in the ring $R$ we get, $I=\langle h_1',h_2',h_3,...,h_d\rangle  $. Define $ \theta=\tau\perp I_{d-2}\in \GL_d(R/I)$. Then note that $(h_1,h_2,h_3,...,h_d)\theta = (h_1',h_2',h_3,...,h_d)$ and $\det(\theta)-f^2\in I$. Since $\det(\theta)-\det(\alpha)\in I$, there exists $\epsilon'\in \SL_d(R/I)$ such that $\theta \epsilon'= \alpha$. Since $\dim(R/I)=1$, by Theorem \ref{MKDLE} the natural map $\SL_d(R)\sur \SL_d(R/I)$ is surjective. Therefore we can lift $\epsilon'$ and get $\epsilon\in \SL_d(R)$ such that they are equal modulo $I$. Let $(G_1,...,G_d)=(h_1',h_2',h_3,...,h_d)\epsilon$. Then note that $I=\langle G_1,...,G_d\rangle  $. It only remains to show $G_i-g_i\in I^2$.
	\par Consider any $d$-tuple $[(a_1,...,a_d)]$ as a map $(R/I)^d\to I/I^2$ sending $e_i\to a_i \mod (I^2)$. Then we have the following. 
	\begin{align*}
		[(G_1,...,G_d)]&=[(h_1',h_2',h_3,...,h_d)\epsilon]\\
		&=[(h_1,h_2,h_3,...,h_d)\theta\epsilon]\\
		&=[(h_1,...,h_d)\theta\epsilon']\\
		&=[(h_1,...,h_d)\alpha]\\
		&=[(f_1,...,f_d)\alpha]\\
		&=[(g_1,...,g_d)]
	\end{align*}This completes the proof.\qed  
	
	\smallskip 
	
	We end this section with a lemma whose proof can be found in \cite[Lemma 5.5]{MPHIL}. The hypothesis in \cite{MPHIL} 
	is $\dim(A)=\hh(J)=n$, but  the same proof works for the following version.

	\bl\label{MKLLM}
	Let $A$ be a Noetherian ring and $J \subset A$ be an ideal of height $n$. Let
	$f \in A\setminus \{0\}$ such that $J_f$ is a proper ideal of $A_f$. Assume that $J_f = (a_1,..., a_n)$, where $a_i \in J$. Then,
	there exists $\sigma \in\SL_n(A_f)$ such that $(a_1,...,a_n)\sigma = (c_1,..., c_n)$, where $c_i=(b_i)_f$ for some $b_i\in J\subset A$ with $ht (\langle b_1,..., b_n\rangle  A) = n$.
	\el
	
	\section{A splitting criterion in terms of generic sections}\label{5}
	
In this section, we prove Theorem \ref{MT1I}, as mentioned in the Introduction. As a preparation, we first deduce a cancellation result for projective modules in dimension two over a certain $C_1$ field with characteristic $\not= 2$. If the ring is smooth, then it is essentially contained in \cite[Theorem 2.4]{AAS}. Here we drop the smoothness assumption of A. A. Suslin's proof following P. Raman's observation (see \cite[Proposition 3.1]{RV}), which is crucial in our set-up.

	\bp \label{SCM}
Let $k$ be an algebraically closed field of characteristic $\not=2$. Let $R$ be an affine algebra over the $C_1$ field $k(T)$ of dimension $2$. Then any projective $R$-module $P$ of rank $2$ is cancellative. In particular, for any $(a,p)\in \Um(R\oplus P)$ there exists a $\sigma\in \Aut(R\oplus P)$ such that $\sigma(a,p)=(1,0)$.
\ep

\proof First, we consider the case when $P$ is free. This case is essentially contained in \cite[Proposition 3.1]{RV}. For the sake of completeness, we provide some details. One may note that since $k$ is algebraically closed, by Tsen’s theorem \cite{Tsen} it follows that $k(T)$ is a $\mathrm{C_1}$ field. Hence cohomological dimension of $k(T)$ is $\le 1$. Since $\text{char}(k(T))=\text{char}(k)\not = 2$ and $\dim(R)=2$, it follows that $R$ satisfies condition (a) in \cite[Proposition 3.1]{RV}. Therefore, applying the same the theorem concludes.

Now we assume that $P$ is an arbitrary projective $R$-module of rank $2$. In this case, the proof essentially follows the argument of \cite[Theorem 4.1, see also Remark 4.2]{BC1}. Hence, we only provide a sketch.

 First, we observe that, there exists an affine $k[T]$-algebra $A$ such that $S^{-1}A=R$, where $S=k[T]\setminus \{0\}.$ Let $\eta$ be the nilradical of $R$. Since the canonical map $\text{Aut}(R\oplus P)\sur \text{Aut}(R/\eta\oplus P/\eta P)$ is surjective, it is enough to assume that $R$ is reduced. In particular, we may assume that $A$ is reduced.  Let $\MI_A$ be the ideal defining the singular locus of $A$. As $k$ is algebraically closed, we have $\hh(\MI_A)\ge 1$, implying that there exists a non-zero divisor $b\in R$ such that $R_b$ is smooth. Further, since $P$ is finitely generated, multiplying $b$ by a suitable non-zero divisor  we may assume that $P_b$ is free. Then there exists $n\in \mathbb N$ such that if we take $s=b^n$, then $sP\subset R^2$.

Let `bar' denote going modulo $\langle s \rangle$. Since $s$ is a non-zero divisor one may observe that it follows from the Bass cancellation that $(\ol a, \ol p)\equiv (\ol 1,\ol 0)\mod \text{E}(R\oplus P)$. Since the canonical map $\text{Aut}(R\oplus P)\sur \text{E}(\ol R\oplus \ol P)$ is surjective we may alter $(a,p)$ suitably by an automorphism and further assume that $p\in sP\subset R^2$ and $a-1\in \langle s \rangle$. If we take $p=(b,c)\in R^2$, then we note that $(a,b,c)\in \Um_3(R)$.  Moreover, using the prime avoidance lemma we may replace $a$ by $a+\lambda_1b+\lambda_2c$ and assume that $a\not\subset \bigcup_{p\in \ms} p$,  where $\ms$ is the set of all minimal prime ideals of $R$. Hence we may further assume that $a$ is a non-zero divisor such that $a-1\in \langle s\rangle$.

Let $C=R/\langle a \rangle $ and let `tilde' denote going modulo $\langle a \rangle$. Now we note that, since $a-1\in \langle s \rangle$ and $sP\subset R^2$, we have $\widetilde{P}=C^2$. Then it follows using the argument given in \cite[Proposition 3.1]{RV} that $\text{SK}_1(C)$ is a $2$-divisible group. Then there exists $\Gamma'\in \SL_2(C)\cap\text{ESp}_4 (C)$ and $v,w\in R$ such that $\Gamma'(b,c)=(v,w^2)$. Then by \cite[Corollary 3.3]{BC1} there exists a $\Gamma\in \text{Aut}(P)$ such that $\widetilde\Gamma=\Gamma'$. Moreover, altering $1\perp \Gamma$ by a multiple of an element in $\text{E}(R\oplus P)$ we may further assume that $\Gamma(a,b,c)=(a,v,w^2)$. Now one may follow the proof of \cite[Theorem 4.1, paragraph 9 onward]{BC1} to conclude that $(a,v,w^2)\equiv (1,0,0)\mod\text{Aut}(R\oplus P)$.\qed

%Applying Theorem \ref{SCM} we obtain that the unimodular row $(a,b,c)$ can be elementarily completed to a unimodular row $(u,v,w^2)$, for some $(u,v,w)\in \Um_3(R)$.

%
%The proof is essentially contained in \cite[Theorem 4.1]{BC1}. One just needs to remove the hypothesis that the field is ``perfect" (see \cite[Remark 4.2]{BC1}). To do this here we argue as follows: first we point out that the perfectness of the field is only used to ensure that there exists a non-zero divisor $s\in R$ such that $R_s$ is smooth. But in our situation since $R$ is an affine algebra over $k(T)$ there exists an affine $k$-algebra $A$ such that $S^{-1}A=R$, where $S=k[T]\setminus \{0\}.$ Let $\eta$ be the nilradical of $R$. Since the canonical map $\Aut(R\oplus P)\surj \Aut(R/\eta\oplus P/\eta P)$ is surjective, it is enough to assume that $R$ is reduced. In particular, this will imply that $A$ is reduced.  Let $\MI_A$ be the ideal defining the singular locus of $A$. As   $k$ is algebraically closed, we have $\hh(\MI_A)\ge 1$.
%
% We choose a non-zero divisor $a\in \MI_A$. Now if we consider the non-zero divisor $t=\frac{a}{1}\in S^{-1}\MI_A$, then $R_t=(S^{-1}A)_{\frac{a}{1}}=S^{-1}(A_a)$. Since localization preserves smoothness, the ring $R_t$ is smooth. We may further alter $t$ with some suitable multiple to make the additional assumption that $P_t$ is a free $R_t$-module of rank $2$. Now one may follow \cite[Theorem 4.1]{BC1} to conclude the proof.
%\qed

%%%%%%%%%%%%%%%%%%%%%%%%%%%%%%%%%%%%%	
	
	\bl\label{WECI}
	Let $k$ be an algebraically closed field of characteristic $\not=2$. Let $R$ be an affine algebra over $k(T)$ of dimension $d$, where $d\ge 2$. Let $I\subset R$ be an ideal such that $\hh(I)=\mu(I)=d$.  Then any set of generators of $I/I^2$  has a lift to a set of generators of $I$.	
	\el
	\proof Let $I=\langle a_1,...,a_d\rangle  $ and assume that $I=\langle f_1,...,f_d\rangle  +I^2$ is given. We prove that there exist $F_1,\cdots,F_d\in I$ such that $I=\langle F_1,...,F_d\rangle  $ where $F_i-f_i\in I^2$, $i=1,\cdots,d$.

	Replacing $(a_1,...,a_d)$ by $(a_1,...,a_d)\epsilon$ for a suitable $\epsilon\in E_d(R)$, we may always assume that $\dim(R/\langle a_3,...,a_d\rangle  )\le 2$.  Let $B=R/\langle a_3,...,a_d\rangle  $. Let `bar' denote reduction modulo $\langle a_3,...,a_d\rangle  $. As an $R/I$-module two sets of generators of $I/I^2$ must differ by some invertible matrix $\alpha\in \GL_d(R/I)$. Let $\det(\alpha)=a\in (R/I)^*$ and $b\in R$ be such that $ab-1\in I$. Then it follows from Proposition \ref{SCM} that the unimodular row $(\ol b, \ol a_2,-\ol a_1)$ is completable to an invertible matrix in $B$.  Applying \cite[Lemma 5.2]{SMBB3} we get $\tau\in M_2(B)$ such that $(\ol a_1,\ol a_2)\tau =(\ol a_1',\ol a_2')$, where $\ol I= \langle \ol a_1',\ol a_2'\rangle  $ and $\det(\tau)-a\in I$.
	
	Coming back to the ring $R$ we gather that $I=\langle a_1',a_2',a_3,...,a_d\rangle  $. We define $ \theta=\tau\perp I_{d-2}\in \GL_d(R/I)$ (here $I_{d-2}$ is the identity matrix of size $d-2$). Then note that $(a_1,a_2,a_3,...,a_d)\theta = (a_1',a_2',a_3,...,a_d)$ and $\det(\theta)-a\in I$. Since $\det(\theta)-\det(\alpha)\in I$, there exists $\epsilon'\in \SL_d(R/I)$ such that $\theta \epsilon'= \alpha$. As $\dim(R/I)=0$, the natural map $E_d(R)\sur \SL_d(R/I)=E_d(R/I)$ is surjective. Therefore, we can lift $\epsilon'$ and get an $\epsilon\in E_d(R)$ such that they are equal modulo $I$. Let $(F_1,...,F_d)=(a_1',a_2',a_3,...,a_d)\epsilon$. Then note that $I=\langle F_1,...,F_d\rangle  $. It only remains to show that $F_i-f_i\in I^2$.
	
	To see this consider any $d$-tuple $[(a_1,...,a_d)]$ as a map $(R/I)^d\to I/I^2$ sending $e_i\to a_i \mod(I^2)$. Then we have $[(F_1,...,F_d)]=[(a_1',a_2',a_3,...,a_d)\epsilon]=[(a_1,a_2,a_3,...,a_d)\theta\epsilon]=[(a_1,...,a_d)\theta\epsilon']=[(a_1,...,a_d)\alpha]=[(f_1,...,f_d)]$. This completes the proof.\qed  
	
	\smallskip
	
	The above proof is inspired by \cite[Lemma 5.3]{SMBB3}. Following their proof and using Proposition \ref{SCM} in place of a cancellation theorem of Murthy-Swan \cite[Theorem 4]{MuSw}, one can easily obtain the following.
	
\bl\label{lift1}
	Let $k$ be an algebraically closed field of characteristic $\not=2$. Let $R$ be an affine algebra over $k(T)$ of dimension $d$, where $d\ge 2$. Let $I\subset R$ be an ideal such that $\hh(I)=d$. Let $L$ be a projective $R$-module of rank one. Assume that there is a surjection $\alpha:L\oplus R^{d-1}\sur I$. Then any surjective map $\ol{\theta}:L\oplus R^{d-1}\sur I/I^2$ can be lifted to a 
	surjection $\phi:L\op R^{d-1}\sur I$ (meaning, $\phi\otimes (R/I)=\ol{\theta}\ot (R/I)$).
	\el
	
	%The next corollary is just a translation of the previous theorem in terms of the literature of ``Euler class group'' defined in \cite[Section 4]{SMBB3}. Here we would like to mention that the hypothesis ``the ring containing $\mathbb{Q}$" in \cite[Section 4]{SMBB3} was not actually required to define the ``Euler class group". It has only been used to assign the ``Euler class'' of a projective module. Hence one can talk about the Euler class group of any ring.
	
	%\bc
	%Let $k$ be an algebraically closed field of $\text{char}(k)\not=2$. Let $R$ be a $d-$dimensional affine algebra over $k(T)$, where $d\ge 2$. Let $I\subset R$ be an ideal such that $\mu(I/I^2)=\hh(I)=d$. Then for any local orientations $\omega_1$ and $\omega_2$ of $I$ we have $(I,\omega_1)=(I,\omega_2)$ in $E^d(R)$.
	%\ec

	We are now ready to prove Theorem \ref{MT1I}, as mentioned in the Introduction.
	
%	 For simplicity, first we give a proof in the case when the projective module has a trivial determinant. Then we show that the general case can be deduced using Proposition \ref{SCM}.

%	
%	\bt\label{BRQRLTd}
%	Let $k$ be an algebraically closed field of $\text{char}(k)\not=2$. Let $R$ be an affine algebra over $k$ of dimension $d$, where $d\ge 2$. Let $P$ be a projective $R[T]$-module of rank $d$ with a trivial determinant. Let $I\subset R[T]$ be an ideal of height $d$ such that there is a surjection $\phi:P\surj I$. If $\mu(I)=d$, then $P$ has a unimodular element. 
%	\et 
%	\proof  \qed 
%	

	\bt\label{BRQRLT}
	Let $k$ be an algebraically closed field of $\text{char}(k)\not=2$. Let $R$ be an affine algebra over $k$ of dimension $d$, where $d\ge 2$. Let $P$ be a projective $R[T]$-module of rank $d$. Let $I\subset R[T]$ be an ideal of height $d$ such that there is a surjection $\phi:P\sur I$. If $\mu(I)=d$, then $P$ has a unimodular element. 
	\et 
	
	\proof We give the proof in two cases. First, we give the proof in the case when $P$ has a trivial determinant. Then we show that the general case can be deduced using Proposition \ref{SCM}.
	
	\paragraph{\textbf{Case - 1}} In this case we us assume that $P$ has a trivial determinant. Note that if $I$ contains a monic polynomial $f\in R[T]$, then $P_f$ has a unimodular element via the map 
	$\phi\otimes R[T]_f$. Hence using ``monic inversion principle" for projective modules, as stated in \cite[Theorem 3.4]{BRS01}, it follows that 
	$P$ has a unimodular element. Therefore,  we  assume that $I$ does not contain any monic polynomial.
	
	We note that, again by \cite[Theorem 3.4]{BRS01}, it is enough to find  a monic polynomial $f\in R[T]$, such that $P_f$ has a unimodular element. We fix a trivialization $\chi:R[T]\cong \wedge^dP$.  Since $k$ is a algebraically closed field, the quotient field of $R$ is $k$. We take $A=S^{-1}R[T]$, where $S=k[T]\setminus \{0\}$. Then one may observe that $A$ is an affine $k(T)$-algebra of dimension $d$. Let $IA$ be the extension of the ideal $I$ in the ring $A$. Therefore, in the ring $A$ we have $\hh(IA)=\mu(IA)=\mu(IA/I^2A)=d$. Applying Lemma \ref{WECI} we get that any set of generators of $IA/I^2A$ can be lifted to a set of generators of $IA$. In particular, the set of generators induced by the triplet $(P\otimes A,\phi\otimes A, \chi\otimes A)$ lifts to a set of generators of $IA$. Therefore, using subtraction principle as stated in \cite[Corollary 3.4]{SMBB3} we conclude that the projective module $S^{-1}P$ has a unimodular element. Since $P$ is finitely generated, we obtain a monic polynomial $f\in k[T]\subset R[T]$, such that $P_f$ has a unimodular element. This concludes the proof when the determinant of $P$ is trivial.

\paragraph{\textbf{Case - 2}} Now we give the proof in the case where $P$ has an arbitrary determinant, say $\det(P)=L$. Let $A=S^{-1}R[T]$, where $S=k[T]\setminus \{0\}$. Following Case - 1, without loss of generality we may assume that assume that $I$ does not contain any monic polynomial. Write $P'= P\otimes A$, and $I'=IA$. Moreover, as observed in the Case - 1, it is enough to show that $P'$ has a unimodular element. 
	
 We fix an isomorphism $\chi:L\simeq \wedge^d(P)$. We write $L'=L\otimes A$ and $\chi'=\chi\otimes A$. The triple $(P',\chi',\phi)$ will induce a surjection $\omega:L'\oplus A^{d-1}\sur I'/I'^{2}$.
	
We then have $\mu(I')=\hh(I')=d$.	Let $I'=\langle a_1,\cdots,a_d\rangle$. Performing elementary transformations we can ensure that 
the ring $B:=A/\langle a_3,\cdots,a_d\rangle$ has dimension $2$ (we also note that $B$ is an affine algebra over the $C_1$ field 
$k(T)$). Let `bar' denote reduction modulo $\langle a_3,\cdots,a_d\rangle$. We can now follow the proof of \cite[Theorem 6.8]{SMBB3}
(third paragraph onward) to conclude that $(\overline{I'})=0$ in $E^2_0(B,\overline{L'})$ (where $E^2_0(B,\overline{L'})$ is the weak Euler class group of $B$). As $\dim(B)=2$, it follows from 
\cite[Proposition 6.2]{SMBB3} that $\overline{I'}$ is a surjective image of a projective $B$-module $Q$ such that $Q$ is stably isomorphic to $\overline{L'}\oplus B$. Then by Proposition \ref{SCM}, $Q$ is isomorphic to   $\overline{L'}\oplus B$. Therefore,
we have a surjection $\overline{\alpha}: \overline{L'}\oplus B\sur \overline{I'}$. Consequently, we get that
$\alpha$ is a surjection from $L'\oplus A^{d-1}\sur I'$. 

Now, applying Lemma \ref{lift1} we conclude that $\omega:L'\oplus A^{d-1}\sur I'/I'^{2}$ has a surjective lift, say,
$\theta: L'\oplus A^{d-1}\sur I'$. We can now apply \cite[Corollary 3.4]{SMBB3} and obtain that $P'$ has a unimodular element. This concludes the proof.
\qed

	\bc\label{EUESFM}
	Let $k$ be an algebraically closed field of $\text{char}(k)\not=2$. Let $R$ be an affine algebra over $k$ of dimension $d$, where $d\ge 2$.  Let $P$ be a stably free $R[T]$-module of rank $d$. Then $P$ has a unimodular element.
	\ec
	\proof  Since $P$ is stably free $R[T]$-module of rank $d$, we have a short exact sequence
	$$
		0 \lra  (R[T])^d \lra   R[T]\oplus P \stackrel{(b,-\alpha)}\lra  R[T]\lra  0,$$
		for some suitably chosen $(b,-\alpha)\in R[T]\oplus P^*$. Using a theorem due to Eisenbud-Evans \cite{EE} (or one can see \cite[Corollary 2.13]{SMBB3}), we may replace $\alpha$ with $\alpha+b\gamma$, for some $\gamma\in P^*$, and assume that either $\hh(\alpha(P))= d$ or $\alpha(P)=R[T]$. Here we note that by this replacement of $\alpha$, the kernel remains unaltered. If $\alpha(P)=R[T]$, then this proves the theorem. So let us assume that $\alpha(P)=I\subset R[T]$ is a proper ideal of height $d$. Then it follows from \cite[Lemma 2.8 (i)]{SMBB3} that $\mu(I)=d$. We can now apply Theorem \ref{BRQRLT} to conclude the proof.\qed
	
	\smallskip

	\section{Monic inversion principle for modules: A question of Roitman}
	
	In this section we give an affirmative answer to a question asked by M. Roitman. The following proposition is crucial to our proof. The proof  is motivated from \cite[3.3, 3.4]{SMBB3}.
	
	\bp\label{SPZ}
	Let $A$ be a finitely generated $\mathbb{Z}$-algebra of dimension $d\ge 3$ such that there exists an integer $n\ge 2$ with $\frac{1}{n}\in A$. Let $P$ be a projective $A$-module of rank $d$ with trivial determinant. Let $\chi :A\cong \wedge^dP $ be an isomorphism. Let $I,J\subset A$ be two comaximal ideals such that $\hh(I)\ge d-1$ and $\hh(J)\ge d$. Moreover, assume that there exist surjections $\alpha:P\sur I\cap J$ and $\beta:A^d\sur I$. Let `bar' denote going modulo $I$. Suppose that, there exists an isomorphism $\delta:\ol A^d\cong \ol P$ with the following properties: 
	\begin{enumerate}
		\item $\ol{\beta}=\ol\alpha\delta$;
		\item $\wedge^d\delta=\ol\chi$.
	\end{enumerate}
	Then, there exists a surjection $\gamma:P\sur J$ such that $\gamma\otimes A/J=\alpha\otimes R/J$.
	\ep
	\proof Let $\beta $ correspond to the set of generators $I=\langle a_1,...,a_d\rangle  $. Here we observe that $\dim(A/J^2)\le 0$. This will imply that any unimodular row in $A/J^2$ of length $d$ can be completed to the first row of an elementary matrix. Since the canonical map $E_d(A)\sur E_d(A/J^2)$ is surjective, going modulo $J^2$ we may assume the following:
	\begin{enumerate}
		\item [(i)] $\langle a_1,...,a_{d-1}\rangle  +J^2=A$;
		\item [(ii)] $a_d\in J^2$.
	\end{enumerate}
	Moreover, in view of the prime avoidance lemma, we may replace $a_i$ by $a_i+\lambda_ia_d$ for some $\lambda_i\in A$, $(i=1,...,d-1)$ and further assume that (iii) $\hh(\langle a_1,...,a_{d-1}\rangle  )=d-1$. Here we note that, we are not changing the notation of $a_i$s'. From (i) it follows that there exists $\lambda\in \langle a_1,...,a_{d-1}\rangle  $ such that $\lambda-1\in J^2$. After replacing $a_d$ by $a_d+\lambda$ we may also assume that (iv) $a_d-1\in J^2$.
	\par Consider the following ideals in $A[T]$: $K'=\langle a_1,...,a_{d-1},a_d+T\rangle  $, $K''=J[T]$ and $K=K'\cap K''$. Then it is enough to show that there exists a surjection $\theta:P[T]\sur K$ such that $\theta(0)=\alpha$. For, if we can achieve this, then specializing at $1-a_{d}$ we get $\gamma:=\theta(1-a_d):P\sur J$. Since $a_d-1\in J^2$, we have $\gamma\otimes A/J=\theta(0)\otimes A/J=\alpha\otimes A/J$. In the rest of the proof we will find such a $\theta$.
	
	As $\dim(A[T]/K')=\dim(A/\langle a_1,...,a_{d-1}\rangle  )\le 1$, the module $P[T]/K'P[T]$ is a free $A[T]/K'$-module of rank $d$. We choose an isomorphism $$\kappa(T):(A[T]/K')^d\cong P[T]/K'P[T]$$ such that $\wedge^d\kappa(T)=\chi\otimes A[T]/K' $. Therefore, we get $\wedge^d\kappa(0)=\wedge^d\delta$. This will imply that $\kappa(0)$ and $\delta$ differ by an element $\alpha'\in \SL_d(A/I)$. It follows from \cite[Theorem 16.4]{SV} that $SK_1(A/I)$ is trivial. As a consequence we get $\SL_d(A/I)=E_d(A/I)$ (as $d\ge 3$). Hence we can lift $\alpha'$ to get some $\alpha\in E_d(A)$, and use this to alter $\kappa(T)$ so that $\kappa(0)=\delta.$
	\par Sending the canonical basis vectors to $a_1,...,a_{d-1},a_d+T$ respectively we can define a surjection from $(A[T])^d\sur K'$. This will induce a surjection $\epsilon(T):(A[T]/K')^d\sur K'/K'^2$. Let $\phi(T):=\epsilon(T)\kappa(T)^{-1}:P[T]/K'P[T]\sur K'/K'^2$. Note that $$\phi(0)=\epsilon(0)\kappa(0)^{-1}=\alpha\otimes A/I.$$Since $d\ge \dim(A[T]/K')+2=3$, using \cite[Theorem 2.3]{MR} we get a surjection $\theta(T):P[T]\sur K$ such that $\theta(0)=\alpha.$ This completes the proof. \qed
	
	\smallskip
	
	%\rmk 
	%
	% The only place at which the hypothesis $d\ge 3$ is used in the proof, is to establish the fact that there exists a natural surjection $\SL_d(A)\surj \SL_d(A/I)$. But this fact can be obtained using the same arguments given in \cite[ Corollary 2.3]{MKD}. Hence the above Lemma can be proved for $d\ge 2$.

	%\bc\label{SPZC}
	%Let $A$ be a finite $\mathbb{Z}-$algebra of dimension $d\ge 3$. Moreover, assume that there exists an integer $n\ge 2$ such that $n\in A^*$. Let $P$ be a projective $A$-module with trivial determinant of rank $d$. Let $\chi :A\cong \wedge^dP $ be an isomorphism. Let $I\subset A$ be an ideal of height $\ge d-1$ such that there exist surjections $\alpha:P\surj I$ and $\beta:A^d\surj I$. Let `bar' denote going modulo $I$. Suppose that, there exists an isomorphism $\delta:\ol A^d\cong \ol P$ such that
	
	%\begin{enumerate}
	%	\item [(i)] $\ol{\beta}=\ol\alpha\delta$;
	%	\item [(ii)]$\wedge^d\delta=\ol\chi$.
	%\end{enumerate}
	%Then, P has a unimodular element.
	%\ec
	%\proof The proof follows from Theorem \ref{SPZ}, taking $J=A$.\qed
	
	\bt\label{RQOZ}
	Let $A$ be a finitely generated $\mathbb{Z}$-algebra of dimension $d\ge 1$ such that there exists an integer $n\ge 2$ with $\frac{1}{n}\in A$.  Let $P$ be a projective $A[T]$-module with trivial determinant of rank $d$. Let $J\subset A[T] $ be an ideal of height $d$ containing a monic polynomial. Assume that there exists a surjection $\alpha:P\sur J$. Then $P$ has a unimodular element.
	\et
	
	\proof For $d=1$ the theorem follows trivially as $\det(P)\cong A$ and for $d=2$ the proof is done in \cite[Proposition 3.3]{BMI}. Therefore, we may assume that $d\ge 3$.  
	
	Fix  an isomorphism $\chi:A[T]\cong \wedge^dP$. Let `bar' denote going modulo $J$. Since $P$ has a trivial determinant and $\dim(A[T]/J)\le 1$, the module $P/IP$ is a free $A[T]/I$-module of rank $d$. Let $\delta:(A[T]/J)^d\cong P/JP$ be an isomorphism such that the following hold.
	\begin{equation}
		\wedge^d\delta=\chi\otimes A[T]/J
	\end{equation}
	We define $\omega:=(\alpha\otimes A[T]/J)\delta:(A[T]/J)^d\sur J/J^2$. Since $J$ contains a monic polynomial and $ \dim(A[T]/J)+2\le 3\le d$, it follows from \cite{NMK} and \cite{SM1} that there exists a surjection $\beta:(A[T])^d\sur J$ such that $\beta \otimes A[T]/J=\omega$.
	
	Since $\wedge^dP$ is extended from the ring $A$, in a view of \cite[Theorem 2.3]{BRS01} it is enough to show that $P/TP$ has a unimodular element. We will establish this with the help of Proposition \ref{SPZ}. Let us define the following notations:
	
	\begin{itemize}
		\item $P(0):=P/TP$, $J(0):=J\otimes A[T]/\langle T\rangle $;
		\item $\alpha(0):=\alpha\otimes A[T]/\langle T\rangle :P(0)\sur J(0)$;
		\item $\beta(0):=\beta\otimes A[T]/\langle T\rangle :A^d\sur J(0)$;
		\item $\omega(0):=\omega\otimes A/J(0):(A/J(0))^d\sur J(0)/J(0)^2$;
		\item $ \delta(0):=\delta\otimes A[T]/\langle T\rangle:(A/J(0))^d\iso P(0)/J(0)P(0)$;
		\item $\chi(0):=\chi\otimes A[T]/\langle T\rangle :A\iso \wedge^d P(0)$.
	\end{itemize}
	Since $J\cap A\subset J(0)$, we have $\hh(J(0))\ge d-1$. Now we observe that
	$\beta(0)\otimes A/J(0)=(\alpha(0)\otimes A/J(0))\delta(0)$, as both matches with $\omega(0)$. It follows from (1) that $\wedge^d\delta(0)=\chi(0)\otimes A/J(0).$ This shows that all the hypotheses of Proposition \ref{SPZ} are satisfied. Now we may apply Proposition \ref{SPZ} (taking $J=A$) to obtain that $P(0)$ has a unimodular element. \qed
	
	\section{Monic inversion principle for ideals}
	
	Here we shall give a direct proof of a \emph{``monic inversion principle"} for ideals in $R[T]$, where $R$ is an affine algebra over $\ol{\mathbb{F}}_p$ of dimension $d\ge 2$. Note that if one assumes that $p> (d-1)!$ or that $R$ is smooth, one can obtain a proof using (\cite{MKD1} or \cite{DTZ2}). We make no such assumption and our method is entirely different.
	
	\bt \label{INJECG} 
	Let $R$ be an affine algebra over $\ol{\mathbb{F}}_p$ of dimension $d\ge 2$. Let $I\subset R[T]$ be an ideal such that $\hh(I)=\mu(I/I^2)=d$. Assume that, $I=\langle f_1,...,f_d\rangle  +I^2$. Suppose that there exist $F_i\in IR(T)$ such that $IR(T)=\langle F_1,...,F_d\rangle  $, with $F_i-f_i\in IR(T)^2$, ($i=1,...,d$). Then there exist $g_i\in I$, such that $I=\langle g_1,...,g_d\rangle  $, with $g_i-f_i\in I^2$.
	\et
	\proof We divide the proof into two cases.\\
	%%%%%%%%%%%%%%%%%%%%%%%%%%%%%%%%%%%%%%%%%%%%%%%%%%%%%%%%%%%%%%%%%%%%%%%%%%%%%%%%%%%%%%%%
	{\textbf{Case 1.}} In this case we assume that $d=2$. Here we note that if the ideal $I$ contains a monic polynomial in $T$, then a lift exists by \cite[Theorem 3.2]{MKD}. Therefore, without loss of generality we may assume that $I$ does not contain a monic polynomial.  Let $I=\langle f_1,f_2\rangle  +I^2$ induce the surjection $\omega_I:(R[T]/I)^2\sur I/I^2$, by sending the canonical basis $e_i$ to the image of $f_i$ in $I/I^2$ $(i=1,2)$. Initial part of our  proof in this case essentially follows \cite[Theorem 7.1]{MKD1}.

	Using a standard patching argument there exists a projective $R[T]$-module $P$ with trivial determinant of rank $2$ and a surjection $\alpha:P\sur I$. Fix an isomorphism $\chi:R[T]\cong \wedge^2 P$. Let $\alpha$ and $\chi$ induce $I=\langle f'_1,f'_2\rangle  +I^2$.
	\par Let `bar' denote going modulo $I^2$. We notice that any two set of generators of $I/I^2$ differs by an invertible matrix in $\GL_2(R[T]/I)$. Hence there exists $\ol{\sigma}\in \GL_2(R[T]/I)$ such that $(\ol f_1,\ol f_2)=(\ol f'_1,\ol f'_2 )\ol\sigma $. Let $\det(\sigma)=\ol f$. Using \cite[ Lemma 2.7 and 2.8]{SMBB3} we get another surjection $\beta:P_1 \sur I$ and an isomorphism $\chi_1: R[T]\iso \wedge^2P_1$, where $P_1$ is projective $R[T]$-module of rank $2$ having trivial determinant. Moreover, the triplet $\{P_1,\chi_1,\beta\}$ satisfy the property that if the set of generators of $I/I^2$
	induced by $\beta$, $\chi_1$ and a fixed basis of $R[T]^2$, say $\{\eta_1,\eta_2\}$, is $\ol h_1, \ol h_2$, then
	$(\ol h_1, \ol h_2) = (\ol f'_1, \ol f'_2)\ol\delta$, where $\ol \delta \in \GL_2(R[T]/I)$ has determinant $\ol f$. Let us define $\ol\gamma:=\ol\sigma^{-1}\ol{\delta}\in \SL_2(R[T]/I)$. Then, we note that $(\ol f_1,\ol f_2)\ol{\gamma}=(\ol h_1,\ol h_2)$. Since by Theorem \ref{MKDLE} the canonical map $\SL_d(R[T])\sur \SL_d(R[T]/I)$ is surjective, we can find a $\gamma\in \SL_2(R[T])$ such that $\gamma$ is a lift of $\ol\gamma$.
	
	Now here we observe that it is enough to lift the set of generators $I=\langle h_1,h_2\rangle  +I^2$ to a set of generators of $I$. As suppose that we obtain such a lift $I=\langle a_1,a_2\rangle  $ of $I=\langle h_1,h_2\rangle  +I^2$. Then we define $(g_1,g_2):=(a_1,a_2)\gamma^{-1}$. This implies that $(g_1-f_1,g_2-f_2)=(a_1-h_1,a_2-h_2)\gamma^{-1}\in I^2\times I^2$. The remaining part of the proof in this case is dedicated to find a lift of $\{h_1,h_2\}$.
	\par Since the set of generators $\{f_1,f_2\}$ and $\{h_1,h_2\}$ differs by a matrix in $\SL_d(R[T]/I)$, it follows from the hypothesis of the theorem that the set of generators $\{h_1,h_2\}$ lifts to a set of generators of $IR(T)$. Therefore, applying \cite[Corollary 3.4 ]{SMBB3} we obtain that the projective module $P_1\otimes R(T)$ has a unimodular element. Since $P_1$ is a projective module of rank $2$ having a trivial determinant, we actually obtain that $P_1\otimes R(T)$ is a free module. Now we may apply ``Affine Horrocks theorem" \cite{Q} and obtain that the module $P_1$ is in fact a free module.
	
	Therefore, we may choose an isomorphism $\eta:(R[T])^2\cong P_1$ such that $\wedge^2\eta=\chi_1$. Let us define $H_i=\beta\eta(\eta_i) $, $i=1,2$. Then note that $I=\langle H_1,H_2\rangle  $ and $\ol H_1\wedge \ol H_2=\ol h_1\wedge \ol h_2$ in $\wedge^2(I/I^2)$. Hence we can find $ \sigma\in \SL_2(R[T]/I)$ such that $(\ol H_1,\ol H_2)\sigma =(\ol h_1,\ol h_2)$. Now using Theorem \ref{MKDLE} there exists $\tau \in \SL_2(R[T])$ such that $\ol \tau=\sigma$. Let $(g_1,g_2)=(H_1,H_2)\tau$. Then $(\ol g_1,\ol g_2)=(\ol H_1,\ol H_2)\ol \tau=(\ol H_1,\ol H_2)\sigma=(\ol h_1,\ol h_2)$. This concludes the proof in the case of dimension two.

	\paragraph{\textbf{Case 2.}} Here we assume that $d\ge 3$. Again as before we may assume that none of the ideals which are going to appear in this proof, contain monic polynomials. There is a monic polynomial $f\in R[T]$ such that $I_f=\langle F_1,...,F_d\rangle  $, with $F_i-f_i\in I^2_f$. Let $B=R[T]/\langle f\rangle  \cap I^2$ and `bar' denote going modulo $\langle f\rangle  \cap I^2$. Note that $\dim(B)\le d$. In  $B$ we have $\ol{I}=\langle \ol f_1,...,\ol f_d\rangle  +\ol I^2$. Using \cite[Theorem 2.4]{MKD} we can find $\ol h_i\in \ol I$ such that $\ol I=\langle \ol h_1,...,\ol h_d\rangle  $, with $\ol f_i-\ol h_i\in \ol I^2$. Therefore,  $I=\langle h_1,...,h_d\rangle  +I^2\cap \langle f\rangle  $. Using Lemma \ref{MKL} there exits $e\in I^2\cap \langle f\rangle  $ such that $I=\langle h_1,...,h_d,e\rangle  $, and $e(1-e)\in \langle h_1,...,h_d\rangle  $. Using prime avoidance lemma (cf. \cite[Corollary 2.13]{SMBB3}) replacing $h_i$ by $h_i+e\lambda_i$ we may assume that $\hh(\langle h_1,...,h_d\rangle  )_e=d$ or $\langle h_1,...,h_d\rangle  _e=R[T]_e$. Note that if $\langle h_1,...,h_d\rangle  _e=R[T]_e$, then $e\in \sqrt{\langle h_1,...,h_d\rangle}  $. A local checking ensures that $I=\langle h_1,...,h_d\rangle$. This will establish the theorem in this case. Hence the only non-trivial situation is when $\hh(\langle h_1,...,h_d\rangle  _e)=d$. From now onward we will make this assumption.
	
	We define $I_1=\langle h_1,...,h_d,1-e\rangle  $. Then note that we have the following:
\begin{enumerate}[\quad \quad (a)]
		\item  $I\cap I_1=\langle h_1,...,h_d\rangle  $;
		\item  $I_1+\langle e\rangle  =I_1+\langle f\rangle  =I_1+I=R[T]$;
		\item $\hh(I_1)=d$.
	\end{enumerate}
	(b) and (c) imply that $\hh(I_1)=d$. We observe that $\hh((I_1)_f)=d$. Since $h_i-f_i\in I^2$, to prove the theorem it is enough to lift $I=\langle h_1,...,h_d\rangle  +I^2$ to a set of generators of $I$. Now $I+I_1=R[T]$. Therefore $I\cap I_1=\langle h_1,...,f_d\rangle  $ will induce $I_1=\langle h_1,...,h_d\rangle  +I_1^2$. In view of subtraction principle \cite[Proposition 7.2]{SBMKD} to prove the theorem, it is enough to lift  $I_1=\langle h_1,...,h_d\rangle  +I_1^2$ to a set of generators of $I_1$. 
	
	Note that in the ring $R[T]_f$, the set of generators $I_f=\langle h_1,...,h_d\rangle  R[T]_f+I_f^2$ lifts to a set of generators of $I_f$. Hence again by subtraction principle \cite[Proposition 7.2]{SBMKD} there exist $l_i\in (I_1)_f$, such that $(I_1)_f=\langle l_1,...,l_d\rangle  $, with $h_i-l_i\in ({I_1})_f^2$. Let $k\ge 1$ be an integer such that $f^{2k}l_i\in I_1$, for all $i$. Since $f$ is unit modulo $I_1$, to find a lift of $I_1=\langle h_1,...,h_d\rangle  +I_1^2$  to a set of generators of $I_1$, by Lemma \ref{NLTP} it is enough to lift $I_1=\langle f^{2k}l_1,...,f^{2k}l_d\rangle  +{I_1}^2$ to a set of generators of $I_1$. Therefore, we may replace $l_i$ with $f^{2k}l_i$ and assume that $l_i\in I_1$. The remaining part of the proof is devoted to find such a lift of $I_1=\langle l_1,...,l_d\rangle  +{I_1}^2$.

	Using Lemma \ref{MKLLM} we get $\epsilon \in \SL_d(R[T]_f)$ such that $(l_1,...,l_d)\epsilon=(l_1',...,l_d')$, where $l_i'\in I_1$ and $\hh(\langle l_1',...,l_d'\rangle  R[T])=d$. Let $\langle l_1',...,l_d'\rangle  R[T]=\bigcap_{i=1}^rq_i\bigcap_{i=r+1}^nq_i$ be the reduced primary decomposition, where $q_i$'s are $p_i$-primary ideal in $R[T]$, such that $f\not\in p_i$ for $i\le r$ and $f\in p_i$ for all $i>r$. Since $(I_1)_f=\langle l_1',...,l_d'\rangle  _f$ is a proper ideal of height $d$, we must have $r\ge 1.$ Let $I_2=\bigcap_{i=r+1}^nq_i$. Therefore, we get $\hh(I_2)\ge d$ and $f\in \sqrt{I_2}$. 
	
	We claim that $I_1=\bigcap_{i=1}^rq_i$. To prove this note that $(I_1)_f=\langle l_1',...,l_d'\rangle  R[T]_f=\bigcap_{i=1}^r{q_i}_f$ implies that $I_1 \subset \bigcap_{i=1}^rq_i$. Assume, if possible, that there exists another $p$-primary ideal $q$ (where $p\not=p_i$, for any $i=1,...,r$) in the reduced primary decomposition of the ideal $I_1$. Since $p_f\not= ({p_i})_f$, for $i=1,...,r$, and $(I_1)_f=\bigcap_{i=1}^r{q_i}_f$ we must have $p_f=R[T]_f$. This implies that $f\in p$. But this is not possible as $I_1+\langle f\rangle  =R[T]$ in particular, $p+\langle f\rangle=R[T]$. Therefore, we get $I_1=\bigcap_{i=1}^rq_i$.
	
	Hence we obtain the following:
	\begin{enumerate}[\quad \quad (1)]
		\item $I_2$ contains a monic polynomial in particular, some power of $f$;
		\item $I_1+ I_2=R[T]$ (as going modulo $I_1$ any power of $f$ is a unit);
		\item  $I_1\cap I_2=\langle l_1',...,l_d'\rangle  R[T]$.
	\end{enumerate}
	Note that $(3)$ gives us for any prime ideal $p\supset I_2$, we must have $\mu((I_2)_p)\le d$. Applying Krull's generalized principal ideal theorem we get $\hh(I_2)=d$. Now, (3) will induce $I_1=\langle l_1',...,l_d'\rangle  +I_1^2$ and $I_2=\langle l_1',...,l_d'\rangle  +I_2^2$. Since $I_2$ contains a monic, using \cite[Proposition 3.2]{MKD} $I_2=\langle l_1',...,l_d'\rangle  +I_2^2$ can be lifted to a set of generators of $I_2$. Therefore, applying the subtraction principle \cite[Proposition 7.2]{SBMKD} we can lift $I_1=\langle l_1',...,l_d'\rangle  +I_1^2$ to a set of generators of $I_1$.

	Since $I_1+\langle f\rangle=R[T]$, we have $R[T]/I_1=(R[T]/I_1)_f$, and hence $\epsilon\in \SL_d(R[T]/I_1)$. As $I_1=\langle l_1',...,l_d'\rangle  +I_1^2$ has a lift to a set of generators of $I_1$ and $(l_1,...,l_d)\epsilon=(l_1',...,l_d'),$ using \cite[Lemma 4.1]{MKDIMRN}  we can lift $I_1=\langle l_1,...,l_d\rangle  +I_1^2$ to a set of generators of $I_1$. Therefore, we can lift $I_1=\langle h_1,...,h_d\rangle  +I_1^2$ to a set of generators of $I_1$. This completes the proof.\qed 
	
	\section{Splitting criterion via an obstruction class in an obstruction group}\label{ECLG}
	
	Let $R$ be an affine algebra over $\k$ of dimension $d\ge 2$ and $P$ be a finitely generated projective $R[T]$-module of rank $d$, with  trivial determinant. Let us fix an isomorphism $\chi:R[T]\simeq \wedge^dP$. The purpose of this section is to define an Euler cycle for the triplet $(P,\lambda,\chi)$, in the group $E^d(R[T])$ and show that the Euler cycle govern the splitting problem for $P$, where $\lambda$ is a generic section of $P$. First we will prove some addition and subtraction principles. For most of the proofs in this section we will frequently move to the ring $R(T)$, prove the results in $R(T)$ then using Theorem \ref{INJECG} we will come back to the ring $R[T]$. Some of the results below were proved for Noetherian ring containing $\mathbb{Q}$ in \cite{MKD1}.
	\subsection*{Addition and subtraction principles}\label{SOASP}
	\bp (Addition principle)\label{ap}
	Let $R$ be an affine algebra over $\k$ of dimension $d\ge 2$. Let  $I, J\subset R[T]$ be two
	co-maximal ideals, each of height $d$. Suppose that $I = (f_1,..., f_d)$ and
	$J = (g_1,..., g_d)$. Then $I \cap J = (h_1,..., h_d)$ where $h_i -f_i \in I^2$ and $h_i-g_i\in I^2$.
	\ep
	\proof  Since $\hh(I)=\hh(J)=d$, in the ring $R(T)$ both the ideals $IR(T)$ and $ JR(T)$ are of height $\ge d$. We note that if one of them is of height $>  d$, then there is nothing to prove. So without loss of generality we may assume that each ideal is of height $d$.
	\par Since $I+J=R[T]$, using the Chinese Remainder Theorem we have $I\cap J/(I\cap J)^2\cong I/I^2\oplus J/J^2$. Hence the given set of generators of $I$ and $J$ will induce a set of generators $a_i$'s of $(I\cap J)/(I\cap J)^2$ such that $a_i-f_i\in I^2$ and $a_i-g_i\in J^2$. Thus to prove the theorem it is enough to find a lift of $I\cap J=\langle a_1,...,a_d\rangle  +(I\cap J)^2$ to a set of generators of $I\cap J$. 
	\par In the ring $R(T)$, we have $\hh(I)=\hh(J)=\hh(I\cap J)=\dim(R(T))=d$. Hence applying addition principle as stated in \cite[Theorem 3.2]{SMBB3} we can find $H_i\in (I\cap J)R(T)$ such that $(I\cap J)R(T)=\langle H_1,...,H_d\rangle  R(T)$, with $H_i-f_i\in IR(T)^2$ and $H_i-g_i\in JR(T)^2$. Now using Theorem \ref{INJECG} we are done. \qed
	
	\bp (Subtraction principle)\label{sp}
	Let $R$ be an affine algebra over $\k$ of dimension $d\ge 2$. Let $I, J\subset R[T]$ be two
	co-maximal ideals, each of height $d$. Suppose that $I = (f_1,..., f_d)$ and
	$I \cap J = (h_1,..., h_d)$ where $h_i -f_i \in I^2$. Then there exists $g_i\in J$ such that $J = (g_1,..., g_d)$ with $h_i-g_i\in I^2$.
	\ep
	\proof The proof uses the same arguments as in Proposition \ref{ap} with slight modification, hence we shall only sketch a proof. As before, without loss of generality we may assume that $\hh(IR(T))=\hh(JR(T))=\hh((I\cap J)R(T))=\dim(R(T))=d$. Since $I+J=R[T]$, we get $J=\langle h_1,...,h_d\rangle  +J^2$. Here we observe that, to prove the theorem it is enough to find a lift of $J=\langle h_1,...,h_d\rangle  +J^2$ to a set of generators of $J$. Applying subtraction principle as stated in \cite[Theorem 3.3]{SMBB3} in the ring $R(T)$ we can find $G_i\in JR(T)$ such that $G_i-h_i\in JR(T)^2$. Then as before we can use Theorem \ref{INJECG} to complete the proof. \qed
	
	%Let $R$ be an affine algebra of dimension $d \ge 2$ over $\k$. Let $I \subset R[T]$ be an ideal of height
	%$d $ such that $I/I^2$
	%is generated by $d$  elements. Two surjections $\alpha, \beta : (R[T]/I)^d\surj 
	%I/I^2$ are said to be related if there exists $\sigma\in  \SL_d(R[T]/I)$ such that $\alpha\sigma=\beta$. This defines
	%an equivalence relation on the set of surjections from $(R[T]/I)^d\surj I/I^2$.
	\subsection*{An obstruction group }
	\bp\label{WDNP}
	Let $R$ be affine algebra over $\k$ of dimension $d\ge 2$ and $I\subset R[T]$ be an ideal of height $d$. Moreover suppose that $\alpha$ and $\beta$ are two surjections from $(R[T]/I)^d\sur I/I^2$
	such
	that there exists $\sigma \in \SL_d(R[T]/I)$ with the property that $\alpha \sigma = \beta$. If
	$\alpha$ can be lifted to a surjection $\theta : (R[T])^d\sur  I$ then so is $\beta$.
	\ep
	\proof Since $\dim(R[T]/I)\le 1$, using Theorem \ref{MKDLE}, we can find $\epsilon\in \SL_d(R[T])$, which lifts $\sigma$. Since $\epsilon\in \SL_d(R[T])$ and $\theta$ is a surjection, it follows that $\theta\epsilon:(R[T])^d\sur I$ is also a surjection. Thus it is only remains to show that $(\theta\epsilon)\otimes (R[T]/I)=\beta$. But this follows from the fact that $\epsilon\otimes R[T]/I=\sigma$ and $\theta\otimes R[T]/I=\alpha.$\qed
	
	\smallskip
	
	Now we proceed to define the $d$-th Euler class group of $R[T]$ where $R$ is an affine algebra of dimension $ d \ge 2$.
	
	\bd
	Let $I\subset R[T]$ be an ideal of height $d$ such that $I/I^2$
	is generated by $d$ elements. Let $\alpha$
	and $\beta$ be two surjections from $(R[T]/I)^d\sur I/I^2$. We say $\alpha$ and $\beta$ are related if there
	exists $\sigma\in \SL_d((R[T]/I)^d$ be such that $\alpha \sigma = \beta$. This defines an
	equivalence relation on the set of surjections from $(R[T]/I)^d\sur I/I^2$. Let $[\alpha]$ denote the
	equivalance class of $\alpha$. If $f_1, ..., f_d$ generate $I/I^2$
	, we obtain a surjection $\alpha : (R[T]/I)^d\sur I/I^2$, sending $e_i$ to $f_i$. We say $[\alpha]$ is given by the set of generators $f_1, ..., f_d$ of $I/I^2$.
	\par Let $G$ be the free Abelian group on the set $B$ of pairs $(I, \omega_J )$, where:
	\begin{enumerate}
		\item $I \subset R[T]$ is an ideal of height $d$;
		\item  $\Spec(R[T]/I)$ is connected;
		\item  $I/I^2$ is generated by $d$ elements; and
		\item  $\omega_I : (R[T]/I)^d\sur I/I^2$
		is an equivalence class of surjections $\alpha : (R[T]/I)^d\sur I/I^2$.
	\end{enumerate}
	
	Let $J \subset R[T]$ be a proper ideal. we get $J_i \subset R[T]$ such that $J = J_1 \cap J_2 \cap ... \cap J_r$, where $J_i$
	’s
	are proper, pairwise co-maximal and $\Spec(R[T]/J_i)$ is connected. We shall say that $J_i$ are the connected
	components of $J$.
	\par Let $K \subset R[T]$ be an ideal of height $d$ such that $K/K^2$
	is generated by $d$ elements. Let $K = \cap K_i$
	be the decomposition of $K$ into its connected components. Then note that for every $i$, 
	$\hh(K_i) = d$. Therefore, applying Chinese Remainder Theorem $K_i/K_i^2$
	is generated by $d$ elements. Let $\omega_K:(R[T]/I)^d\sur K/K^2$ be a surjection. Then in a natural way $\omega_K$ gives rise to surjections
	$\omega_{K_i}: (R[T]/K_i)^d \sur K_i/{K_i}^2$. We associate the pair $(K,\omega_K)$, to the element $\sum_{i=1}^{r}(K_i, \omega_{K_i})$ of $G$. We will call $(K,\omega_K)$ as a \textit{local orientation} of $K$ induced by the surjection $\omega_K$.
	\par Let $H$ be the subgroup of $G$ generated by the set $S$ of pairs $(J, \omega_J )$, where $\omega_J:(R[T]/J)^d\sur J/J^2$ has a surjective lift $\theta:(R[T]/J)^d\sur J$. That is, $\theta$ is a surjection such that $\theta\otimes R[T]/J=\omega_J$. Then, we define the quotient group $G/H$ as the $d$-th Euler class group of $R[T]$ denoted as $E^d(R[T])$. A local orientation $(J,\omega_J)$ is said to be a \textit{global orientation} if $\omega_J$ lifts to a set of generators of $J$. 
	
	\smallskip
	
	\ed
%	\rmk Note that the decomposition of $K$ into its connected components is unique by \cite[Lemma 4.5]{MKD1} as the proof of Lemma 4.5 does not require the assumption that the ring contains $\mathbb Q$.
	\smallskip
	
	\rmk The equivalence class of the pair $(I,\omega_I)\in E^d(R[T])$ defined above is well-defined by Proposition \ref{WDNP}.

	\bl\label{LL}
	Let $R$ be an affine $\k$-algebra of dimension $d\ge 2$, and let $I \subset R[T]$ be an ideal of height $d$ such
	that $I/I^2$
	is generated by $d$ elements. Let  $\omega_I : (R[T]/I)^d\sur I/I^2$
	be a local
	orientation of $I$. Suppose that, the image of $(I, \omega_I )$ is zero in the $d$-th Euler class group
	$E^d(R[T])$. Then $I$ is generated by $d$ elements and $\omega_I$ can be lifted to a
	surjection $\theta : (R[T])^d\sur  I$.
	\el
	
	\proof The proof is done in \cite[Theorem 4.7]{MKDIMRN}, with some additional assumption $p\not=2$, whenever $d=3$, but this can be removed if one uses Proposition \ref{sp} in the appropriate places. Hence we opt to skip the proof to avoid repeating similar arguments. \qed
	
	%\bt
	%Let $A$ be a $d-$dimensional affine domain over a field $k$. Moreover assume that one of the following holds:
	%
	%\begin{enumerate}
	%	\item  $d\ge 3$, $k=\ol k$ and either $char(k)=0$ or $char(k)\rangle  d$;
	%	\item $A$ is a real affine domain such that: either there are no real maximal ideals, or the intersection of all real maximal ideals has height at least one and $d\ge 3$;
	%	\item $d\ge 2$ and $k=\k$ ($p\not=2$) (In this case we can drop the assumption ``domain").
	%\end{enumerate}  
	%Then $E^d(A[T])\cong E^d(A(T))$.
	%
	%\et
	
	\subsection*{An obstruction class}\label{OC}
	\smallskip
	
	\bd\label{ECP}
	Let $R$ be an affine algebra over $\k$ of dimension $d\ge 2$, and let $P$ be a projective $R[T]$-module of rank $d$ having trivial determinant.
	Let $\chi : R[T] \cong \wedge^dP$ be an isomorphism. To the pair $(P, \chi)$, we associate an
	element $e(P, \chi)$ of $E^d(R[T])$ as follows:	let $\lambda : P \sur I$ be a surjection, where $I$ is an ideal of $R[T]$ of height
	$d$. We obtain an induce surjection $\lambda\otimes R[T]/I : P/IP \sur  I/I^2$. Note that, since $P$ has trivial determinant and
	$\dim (R[T]/I) \le 1$, the module $P/IP$ is in fact a free $R[T]/I$-module of rank $d$. We
	choose an isomorphism $\phi : (R[T]/I)^d\iso  P/IP$, such that $\wedge^d\phi=\chi\otimes R[T]/I$. Let $\omega_I$ be the surjection $ (\lambda\otimes R[T]/I)\circ\phi: (R[T]/I)^d\sur I/I^2$. We say that $(I, \omega_I)$ is an \textit{Euler cycle} induced by the triplet
	$(P,\lambda,\chi)$. 
	\par Whenever the class of $(I,\omega_I) \in E^d(R[T])$, induced by the triplet $(P,\lambda,\chi)$ becomes independent of a certain choice of the pair $(\lambda,I)$, we shall call the image of $(I,\omega_I) \in  E^d(R[T])$ as the \textit{Euler class} $e(P,\chi)$ of the pair $(P,\chi)$.
	\smallskip
	
	\ed
	\notation Continuing with the above notations we might omit $P$ sometimes and only say $(I,\omega_I)$ is induced by $(\lambda,\chi)$, if there are no confusions.  By saying an Euler cycle induced by $(P,\chi)$ we mean to say an Euler cycle induced by the triplet $(P,\lambda, \chi)$, for some suitably chosen $\lambda\in P^*$. 
	
	In the next theorem we shall show that the vanishing of any Euler cycle of the triplet $(P,\lambda,\chi)$ is necessary and sufficient for the projective module $P$ to have a unimodular element. 
	\bt\label{EUE}
	Let $R$ be an affine $\k$-algebra of dimension $d\ge 2$ and $P$ be a projective $R[T]$-module with trivial determinant of rank $d$. Let $(P,\lambda,\chi)$ induce an Euler cycle $(I,\omega_I)$, where $\lambda$ is a generic section of $P$. Then, the projective module $P$ has a unimodular element if and only if $(I,\omega_I)=0$ in $E^d(R[T])$. 
	\et
	\proof  
	
	Suppose assume that $(I,\omega_I)=0$. Recall that $(P,\lambda,\chi)$ induces $(I,\omega_I)$ means there exists an isomorphism $\phi : (R[T]/I)^d\iso P/IP$, such that $\wedge^d\phi=\chi\otimes R[T]/I$ and $\omega_I=(\lambda\otimes R[T]/I)\circ \phi$.  Note that, in a view of monic inversion principle for modules as stated in \cite[Theorem 3.4]{BRS01} to show that $P$ has a unimodular element it is enough to show that $P\otimes R(T)$ has a unimodular element. The proof is devoted to show this only.
	
	Since $(I,\omega_I)$ vanishes in $E^d(R[T])$, applying Lemma \ref{LL} we can find a surjection $\theta:(R[T])^d\sur I$, such that $\theta\otimes R[T]/I=\omega_I.$ Therefore, in the ring $R(T)$ we have the following:
	\begin{enumerate}[\quad\quad (1)]
		\item $\lambda\otimes R(T):P\otimes R(T)\sur IR(T)$;
		\item  $\theta\otimes R(T):(R(T))^d\sur IR(T)$; and
		\item $\phi : (R[T]/I)^d\iso  P/IP$ 
	\end{enumerate}
	such that $\omega_I\otimes R(T)/IR(T)=(\lambda\otimes R(T)/IR(T))\circ (\phi\otimes R(T)/IR(T))$ and $\wedge^d(\phi\otimes R(T)/IR(T))=\chi\otimes R(T)/IR(T)$. Hence using a subtraction principle as stated in \cite[Corollary 3.4]{SMBB3} it follows that $P\otimes R(T)$ has a unimodular element.
	\smallskip
	
	Conversely, we assume that, $P=Q\oplus R[T]$. Let $\lambda=(\theta, a)\in Q^*\oplus R[T]$. Note that the since $Q$ has trivial determinant, without loss of generality
	we may assume that $\chi$ is induced by an isomorphism $\chi':R[T]\cong \wedge^{d-1}Q$.
	
	Let us denote the image ideal $\theta(Q)=J$. Note that we can always replace $\alpha$ by an automorphism. Using a theorem due to Eisenbud-Evans \cite{EE}, after performing an elementary automorphism we may assume that $\hh(J)\ge d-1$. Since $\dim(R[T]/J)\le 2$, and $Q/JQ$ has trivial determinant, by \cite[ Corollary 2.10]{MKD} there exists an isomorphism $\gamma:(R[T]/J)^{d-1}\cong Q/JQ$ such that $\wedge^{d-1}\gamma'\cong \chi'\otimes R[T]/J$. Let $\Gamma:R[T]^{d-1}\sur J'\subset J$ be a lift of $(\theta\otimes R[T]/J)\circ \gamma$. Then we have $J'+J^2=J$. Applying Lemma \ref{MKL} we can find an $e\in J^2$ such that $J=J'+\langle e\rangle  $. Moreover, since $I=J+\langle a\rangle$, we may further have $ I =J'+\langle b\rangle  $, where $b=e+(1-e)a$. This will produce a surjection $(\Gamma, b):R[T]^d\sur I$. 
	
	The only remaining part is to show $(\Gamma,b)\otimes R[T]/I= \omega_I$, as then it will imply that  $(I,\omega_I)=0$. To show this, observe that we have a canonical isomorphism $\phi:=(\gamma,1):(R[T]/I)^{d}\cong Q/IQ \oplus (R[T]/I)$ such that $\wedge^d\phi=\wedge^{d-1}\gamma\otimes (R[T]/I)=\chi'\otimes (R[T]/I)=\chi\otimes (R[T]/I)$. Since $b-a\in I^2$ and $\Gamma$ is a lift of  $(\theta\otimes R[T]/J)\circ \gamma$, it follows that $(\Gamma,b)$ is a lift of $\omega_I$. This completes the proof. \qed
	
	The same proof will give us the following corollary and we therefore omit the proof.
	\bc
	Let $R$ be an affine algebra over $\k$ of dimension $ d \ge 2$. Let $P$ and
	$Q$ be projective $R[T]$-modules of rank $d$ and $d-1$ respectively, such that their
	determinants are trivial. Let $\chi :\wedge^dP
	\cong \wedge^d(Q \oplus R[T])$ be an isomorphism. Let $I\subset R[T]$ be an ideal of height $d$ such that there exist surjections $\alpha:P\sur I$ and $\beta:Q\oplus R[T]\sur I$. Let `bar' denote going modulo $I$. Suppose that, there exists an isomorphism $\delta:\ol P\cong \ol{Q\oplus R[T]}$ with the following properties:
	\begin{enumerate}[\quad \quad (i)]
		\item $\ol{\beta}\delta=\ol\alpha$;
		\item  $\wedge^d\delta=\ol\chi$.
	\end{enumerate}
	Then $P$ has a unimodular element.
	\ec
	
	\rmk Here we would like to remark that in this section we have been able to prove results,  till now, without the hypothesis 
	``$\frac{1}{(d-1)!}\in R$''. 
	
	\bt \label{WDN}
	Continuing with the notations as in Definition \ref{ECP}, furthermore, assume that $\frac{1}{(d-1)!}\in R$. Then, the assignment sending the pair $(P, \chi)$ to the element e$(P, \chi)\in E^d(R[T])$, as
	described in Definition \ref{ECP}, is well defined.
	\et
	\proof Let $\mu : P\sur  J$ be another surjection where $J \subset R[T]$ an ideal of
	height $d$. Let $(J, \omega_J
	)$ be obtain from $(\mu, \chi)$. Then to prove the theorem we need to show that $(I,\omega_I)=(J,\omega_J)$ in $E^d(R[T])$.
	
	Applying Lemma \ref{ML}, we get an ideal $K\subset R[T]$ such that $K$ is co-maximal with $I$, $J$ and there exists a surjection $\nu:(R[T])^d\sur I\cap K$ such that $\nu\otimes R[T]/I=\omega_I$. Since $I$ and $K$ are co-maximal $\nu$ induces a local orientation $\omega_K$ of $K$. Therefore, we obtain $(I,\omega_I)+(K,\omega_K)=0$ in $E^d(R[T])$.
	
	Let $L=K\cap J$. Then, again as before since $K$ and $J$ are co-maximal $\omega_K$ and $\omega_J$ will induce a local orientation $\omega_L$ of $L$. This will give us $(L,\omega_L)=(K,\omega_K)+(J,\omega_J)$ in $E^d(R[T])$. Thus to prove the theorem it is enough to show that $(L,\omega_L)=0$. Moreover, in view of Theorem \ref{INJECG} it is enough to show that $(L\otimes R(T), \omega_L\otimes R(T))=0$ in $E^d(R(T))$. But since $\frac{1}{(d-1)!}\in R$, by \cite[Section 4]{SMBB3} the Euler class $e(P\otimes R(T),\chi\otimes R(T))$ is well-defined in $E^d(R(T))$. This concludes the proof. \qed

	\smallskip
	
	\rmk Let $R$ be an affine $\k$-algebra of dimension $d\ge 2$. Moreover, assume that $\frac{1}{(d-1)!}\in R$.  Then Theorem \ref{WDN} and Theorem \ref{EUE} say that the assignments of $(P,\lambda,\chi)$ to an element $e(P,\chi)\in E^d(R[T])$, is precisely the obstruction class for $P$ to have a unimodular element. In this case we shall call the Euler class group $E^d(R[T])$ is the obstruction group to detect the existence of a unimodular element in a projective $R[T]$-module (with trivial determinant) of rank equal to the dimension of $R$. In \cite{MKDIMRN}, taking $R$ to be \emph{smooth} and $p\not=2$, M. K. Das defined the $(d-1)$-th Euler class group $E^{d-1}(R)$ and showed that $E^{d-1}(R)$ is the precise  obstruction group to detect the splitting problem of projective $R$-module of rank $d-1$ (with trivial determinant). Therefore, the methods or results of  \cite{MKDIMRN} cannot be applied here.
		
		\smallskip

	In the language of Euler class theory, we have actually proved the following:

	\bc\label{ECGOG}
	Let $R$ be an affine $\k$-algebra of dimension $d\ge 2$ and  $\frac{1}{(d-1)!}\in R$. Let $P$ be a projective $R[T]$-module with trivial determinant (via $\chi$) of rank $d$. Then $P$ has a unimodular element if and only if $e(P,\chi)=0$ in $E^d(R[T])$.
	\ec

	\section{An analogue of Mohan Kumar's result}
	
	\subsection*{On polynomial $\k$-algebras}
	This section is a sequel of Section \ref{5}. Here we have shown that one can strengthen the results of section \ref{5} whenever the base field is $\k$. We begin with the following remark.
	\smallskip 
	
	\rmk Let $R$ be an affine algebra over $\ol {\mathbb F}_p$ of dimension $d\ge 2$. Note that following the arguments given in \cite[Section 6]{MKD1} and using Theorem \ref{INJECG} one can define the $d$-th weak Euler class group $E^d_0(R[T])$ without the assumption that ring contains the field of rationals.

	\bl\label{WECECG}
	Let $R$ be an affine algebra over $\ol {\mathbb F}_p$ of dimension $d\ge 2$, where $p\not= 2$. Let $I\subset R[T]$ be an ideal such that $\hh(I)=\mu(I/I^2)=\mu(I)=d$. Then any set of generators of $I=\langle f_1,...,f_d\rangle+I^2$ lifts to a set of generators of $I$. In particular, $E^d(R[T])\cong E^d_0(R[T])$.
	\el
	\proof Let $S=\k[T]\setminus \{0\}$ and $A=S^{-1}R[T]$. Then, we observe that $A$ is an affine algebra of dimension $d$ over the field $\k(T)$. If $I$ contains a monic polynomial in $R[T]$, then it follows from \cite[Theorem 3.2]{MKD} that $I=\langle f_1,...,f_d\rangle  +I^2$ has a lift. Therefore, without loss of generality we may assume that $I$ does not contain a monic polynomial in $R[T]$. In particular, this will imply that $\hh(IA)=\mu(IA)=d$. Hence applying Lemma \ref{WECI} we get $F_i\in IA$ such that $IA=\langle F_1,...,F_d\rangle  $, with $f_i-F_i\in IA^2$. Since  $I$ does not contain a monic polynomial in $R[T]$, we also have $\hh(IR(T))=\mu(IR(T))=\dim(R(T))=d$. Therefore, applying Theorem \ref{INJECG} we get a lift of $f_i$'s to a set of generators of $I$. This completes the proof.\qed
	
	\bt\label{MKLIPR}
Let $R$ be an affine algebra over $\ol {\mathbb F}_p$ of dimension $d\ge 2$.  where $p\not= 2$. Let $P$ be a projective $R[T]$-module with trivial determinant of rank $d$.  Assume that there is a surjection $\phi:P\sur I$, where  $I\subset R[T]$ is an ideal of height $d$. 
\begin{enumerate}
\item
If $p\not=2$ and $\mu(I)=d$, then $P$ has a unimodular element.
\item
If $P$ has a unimodular element then $\mu(I)=d$.
\end{enumerate}
\et
	
	\proof Since $P$ has a trivial determinant, we fix an isomorphism $\chi:R[T]\cong \wedge^dP$. Let $(I,\omega_I)$ be the Euler cycle induced by the triplet $(P,\chi,\phi)$.
	
	Assume that $\mu(I)=d$. First we would like to mention that this part is essentially contained in Theorem \ref{BRQRLT} (where we do not need $\det(P)$ to be trivial). Here  we give a different proof using the machinery developed in Section \ref{SOASP}. Applying Lemma \ref{WECECG} we get that the canonical homomorphism $E^d(R[T])\to E^d_0(R[T])$ is in fact an isomorphism. Since $\mu(I)=d$, in particular, this will imply that $(I,\omega_I)=0$ in $E^d(R[T])$. Now we may apply Theorem \ref{EUE} to obtain that $P$ has a unimodular element. This proves (1).
	
	Conversely, assume that $P$ has a unimodular element. Therefore, using Theorem \ref{EUE} it follows that $(I,\omega_I)=0$ in $E^d(R[T])$. Hence by Theorem \ref{LL} we obtain that $\omega_I$ is a global orientation of $I$. In particular, this will imply that $\mu(I)=d$. This completes the proof. \qed
	
	\subsection*{On affine $\k$-algebras}\label{aafpbar}
	In the remaining part of the section we show that the ``smoothness" assumption of \cite[Theorem 5.6]{MKDIMRN} can be removed. Before the main result we need some preparation.
	
	%\notation Let $R$ be a ring. For any two matrix $M\in M_{m\times m}(R)$ and $N\in M_{n\times n}(R)$, by $M\perp N$ we denote the matrix $$\begin{pmatrix}
		%	M & {0} \\
		%	0 & N
		%\end{pmatrix}\in M_{{(m+n)}\times {(m+n)}}(R).$$
		
		\bp\label{RSL}
		Let $A$ be a commutative Noetherian ring and $I\subset A$ be an ideal such that $I=\langle a_1,...,a_n\rangle  $. Let $u\in A$ such that the canonical image of $u$ in the ring $A/I$ is a unit. Let $v\in A$ be such that $uv-1\in I$. Let $r$ be an even integer such that $2\le r\le n$. Let $B=A/\langle a_{r+1},...,a_n\rangle  $ and let `bar' denote going modulo $\langle a_{r+1},...,a_n\rangle  $. Furthermore, assume that the unimodular row $(\ol v, \ol a_2, -\ol a_1,...,\ol a_r, -\ol a_{r-1})$ can be completed to an invertible matrix in $B$. Then, there exists $\alpha\in M_{n}(A)$ such that 
		\begin{enumerate}
			\item $\det(\alpha)- u\in I$;
			\item If we define $( b_1,..., b_n):=( a_1,..., a_n)\alpha$, then $I=\langle b_1,...,b_n\rangle  .$
		\end{enumerate}
		\ep
		\proof Let $e_i$ be the canonical basis of $B^{r+1}$. Let $\delta\in\SL_{r+1}(B)$ such that 
		$$ e_1\delta=(\ol v, \ol a_2, -\ol a_1,...,\ol a_r, -\ol a_{r-1}) .$$ Since $\delta\in \SL_{r+1}(B)$, the rows of $\delta$, say $\{e_1\delta,e_2\delta,...,e_{r+1}\delta\}$, form a basis of $B^{r+1}$. Let us define a $B$-linear surjection $$f:B^{r+1}\sur \ol I \text{ by } f(e_1)=\ol 0, \text{ and } f(e_i)=\ol a_{i-1} \text{ for } i\ge 2.$$ Then we note that $f(e_1\delta)=0$. We define $ c_{i-1}:=f(e_i\delta)$, for all $2\le i\le r+1$. This implies, \begin{equation}\tag{A}\ol I=\langle c_1,...,c_r\rangle.\end{equation} By our choice of the matrix $\delta$, it is of the form  
		\begin{equation}\tag{B}
		\delta=\begin{pmatrix}
			\ol v & \ol a_2, -\ol a_1,...,\ol a_r, -\ol a_{r-1} \\
			* & \alpha'
		\end{pmatrix},
	\end{equation}
		where $*\in B^r$ is a column vector and $\alpha'=(\lambda_{ij})\in M_{r}(B)$. Hence we observe that $$(\widehat{e}_i\alpha')(\ol a_1,..., \ol a_r)^T=\sum_{i=1}^r\lambda_{ij}\ol a_i=f(e_{j+1}\delta) \,(\text{ as } f(e_1)=0 \,) =c_i,$$ where $\widehat{e}_i $ is the canonical basis (row) vector of $B^r$, for all $i=1,...,r$. In particular, this implies $(\ol a_1,...,\ol a_r)(\alpha')^T=(c_1,...,c_r)$. We choose $\beta\in M_{r}(A)$ such that $\ol \beta=(\alpha')^T$. Let `tilde' denote going modulo $I$. Since $\delta\in\SL_{r+1}(B) $ we have $\widetilde{\delta}\in\SL_{r+1}(A/I)$. In particular, from (B) it follows that $\det(\beta)-u\in I$.
		
		We define $\alpha:=(\beta\perp \text{I}_{n-r})$, where $\text{I}_{n-r}$ is the identity matrix in $A$ of order $n-r$. Let $(b_1,...,b_n)=(a_1,...,a_n)\alpha$. Here we observe that $b_{r+i}=a_{r+i}$ for all $i=1,...,n-r$. Moreover, from the definition of $\alpha$ we get $\widetilde{\det(\alpha)}=\widetilde{\det(\beta)}=\widetilde u$. Therefore, it only remains to establish (2). To see this we notice that it is enough to show that $\ol I=\langle \ol b_1,...,\ol b_r\rangle$. But since $\beta$ is a lift of $(\alpha')^T$ we have $c_i=\ol b_i$ for $i=1,...,r$. Therefore, from (A) it follows that $\ol I=\langle \ol b_1,...,\ol b_r\rangle$. This completes the proof. \qed
		
		%Then observe that $I=\langle ...,1,...,b_r,a_{r+1},...,a_n\rangle  $. We choose $\delta\in M_{(n+1)\times (n+1)}(A)$ such that $\ol{\delta}=\delta'$. Notice that there exists $\epsilon \in \text{E}_{n+1}(A)$ such that $( v,  a_1,..., a_n)(\delta^T\perp {I}_{n-r})\epsilon=( 0, b_1,..., b_r, a_{r+1},...,a_n)$, where $I_{n-r}$ is the identity matrix in $A$ of order $n-r$. Let 
		%where $C\in A^n$, $\alpha\in M_{n\times n}(A)$ and $(w_0,w)^T\in A^{n+1}$ such that 
		%(\ol w_0 ,\ol w)^T=(\ol v, \ol a_2, -\ol a_1,...,\ol a_r, -\ol a_{r-1},\ol 0,...,\ol 0). 
		%
		%
		% Here we observe that if we can show that $( a_1,..., a_n)\alpha=( b_1,..., b_r,a_{r+1},...,a_n)$, then by taking $b_{r+i}=a_{r+i}$ the prove completes. A further deduction is the following:
		%
		%Suppose we have shown that $(\ol{a}_1,...,\ol {a}_r)\ol \sigma^T=(\ol{b}_1,...,\ol {b}_r)$. Then, we may replace $\alpha$ by $\alpha\epsilon$ for some suitably chosen $\epsilon \in E_n(A)$ and assume that $( a_1,..., a_n)\alpha=( b_1,..., b_r,a_{r+1},...,a_n)$. Note that this alteration of $\alpha$ does not change $\det(\alpha)-u\in I$. Hence in the remaining part of the proof we will show this only.
		%
		%
		%To see this we observe that for all $i=2,...,r+1$ we have $(\ol v,\ol a_1,...,\ol a_r)\delta^Te_i^T=(\ol v,\ol a_1,...,\ol a_r)(e_i\delta)^T=\ol b_i$. $(a_1,..., a_r)\delta^T=( b_1,..., b_r)$ that is to show that $b_i=(a_1,..., a_r)\delta^Te_i^T$ for all $i$. But note that  Hence this complete the proof.\qed

		\bl\label{EEWECG}
		Let $R$ be an affine algebra of dimension $d\ge 5$ over $\k$, where $p\not= 2,3$.
		Let $J \subset R$ be an ideal of height $d-1$ such that $\mu(J) = d-1$. Assume that, there exist $b_i\in J$ ($i=1,...,d-1$) such that
		$J = \langle b_1, ...,b_{d-1}\rangle   + J^2$. Then, there exist $c_1, ...,c_{d-1} \in J$ such that $J = \langle c_1,..., c_{d-1}\rangle  $
		and $b_i - c_i \in J^2$
		for all $ i = 1, ..., d - 1$.
		\el
		\proof Let $J =\langle a_1, ..., a_{d-1}\rangle  $. We consider any $(d-1)$-tuple $[(a_1,...,a_{d-1})]$ as a map $(R/J)^{d-1}\to J/J^2$ sending $e_i\to a_i \mod(J^2)$. Applying \cite[Lemma 2.2]{BCTD2} we can find a matrix $\delta \in \GL_{d-1}(A/J)$ such that
		$[(a_1,... , a_{d-1})\delta] = [(b_1, ..., b_{d-1})]$. Let $\beta\in M_{(d-1)\times (d-1)}(A)$ be a lift of $\delta$. Let $u \in A$ such that $u\det(\beta)-1\in J$. Applying prime avoidance lemma (cf. \cite[Corollary 2.13]{SMBB3}), replacing $a_i$'s suitably we may assume that $B = R/\langle a_5, ..., a_{d-1}\rangle  $ is an affine $\k$-algebra of dimension $5$. Let `bar' denote going modulo $\langle a_5, ..., a_{d-1}\rangle  $. Since $p\ge 5$, using \cite[Theorem 1.2]{MKAD} we get the unimodular row $(\ol{u}, \ol a_2,-\ol a_1, \ol a_4,-\ol a_3) \in \Um_5(B)$ is completable to an invertible matrix in $B$. We may now apply Proposition \ref{RSL} on the unimodular row $(\ol{u}, \ol a_2,-\ol a_1, \ol a_4,-\ol a_3)$ to obtain the following.
		\begin{enumerate}
			\item $\theta \in M_{d-1}(A)$ such that $\det(\theta)-\det(\beta)\in J$;
			\item if we take $(\beta_1, ..., \beta_{d-1}) = (a_1,..., a_{d-1})\theta$, then
			$J = \langle \beta_1,..., \beta_{d-1}\rangle  $.
		\end{enumerate}
		Let `tilde' denote going modulo $J$. Then, from (1) it follows that $(\widetilde{\theta})^{-1}\widetilde{\beta} \in \SL_{d-1}(A/J)$.  We can
		now apply \cite[ Corollary 2.3]{MKD} and lift $(\widetilde{\theta})^{-1}\widetilde{\beta}$ to a matrix $\gamma \in \SL_{d-1}(A)$. We define $(\beta_1,..., \beta_{d-1})\gamma =
		(c_1,..., c_{d-1})$. Then $J = \langle c_1,..., c_{d-1}\rangle  $. It only remains to show that $b_i - c_i \in J^2$
		for all $ i = 1, ..., d - 1$. However, this follows from the following fact.
		\begin{align*}
			[( {c_1},..., {c_{d-1}})]&=[({\beta_1},..., {\beta_{d-1}})\widetilde{\gamma}]\\
		&=[({\beta_1},..., {\beta_{d-1}})(\widetilde{\theta})^{-1}\widetilde{\beta}]\\
		&=[(  {a_1},..., {a_{d-1}})\widetilde\delta]\\
		&=[({b_1},...,{b_{d-1}})]
		\end{align*}This concludes the proof.\qed

		%\rmk For $d=3$ we can always prove this as we can assume that $u$ is a perfect square. Hence using Swan-Towber we are done in this case. Only remaining case is $d=4$.
		\smallskip
		
		\bt
		Let $R$ be an affine algebra of dimension $d\ge 5$ over $\k$, where $p\not= 2,3$. Let $P$
		be a projective $R$-module of rank $d-1$ with trivial determinant. Then $P$ splits off a free summand of rank one if and only if there exists a generic section $I$ of $P$ such that $I$ is generated by $d-1$ elements.  
		\et
		\proof If $P$ has a unimodular element, then the proof is essentially contained in \cite[Theorem 5.1]{MKDIMRN}. Therefore, we will only need to prove the `only if' part. Suppose that, $\alpha:P\sur I$ is a surjection where $\hh(I)=d=\mu(I)$. We fix an isomorphism $\chi:R\cong \wedge^{d-1}P$. Since determinant of $P$ is trivial, it follows from \cite[Corollary 2.10]{MKD} that there exists $\delta:(R/I)^{d-1}\cong P/IP$ such that $\wedge^{d-1}\omega=\chi\otimes R/I$. Let $\omega: (R/I)^{d-1}\sur I/I^2$ be the surjection defined by $\omega=(\alpha\otimes R/I)\circ \delta $. Since $\mu(I)=d-1$, applying Lemma \ref{EEWECG} we observe that any set of generators of $I/I^2$ lifts to a set of generators of $I$. In particular, the map $\omega$ has a surjective lift. Then, using a subtraction principle as stated in \cite[Theorem 3.4]{MKDIMRN} it follows that $P$ has a unimodular element. \qed
		\smallskip
		
		\rmk For $d=2$, the above theorem is true trivially. Therefore, the only remaining cases are $d=3$ and $4$.
	
	%\bibliographystyle{abbrv}
	%\bibliography{Thesisbib} 

	\bibliographystyle{alphaurl}

\end{document}